\documentclass[11pt]{article}
\usepackage{amsmath,amsthm,amsfonts,amssymb,mathrsfs,bm}
\usepackage[usenames]{color}
\usepackage{hyperref}

\parskip 	\bigskipamount
\newtheorem{theorem}{Theorem}[section]

\newtheorem{corollary}[theorem]{Corollary}
\newtheorem{proposition}[theorem]{Proposition}

\newtheorem{remark}{Remark}[section]

\newtheorem{definition}{Definition}

\def\D{\mathcal{D}}
\def\Z1{\m}
\def\Z{\mathbb{Z}}

\def\R{\mathbb{R}}
\def\C{\mathbb{C}}
\def\P{\mathbb{P}}
\def\eps{\varepsilon}
\def\O{\Omega}
\def\A{\mathcal{O}}

\def\b{\widetilde{B}}

\def \del{\delta}
\def \z{\zeta}
\def \o{\omega}
\def \Z{\boldsymbol{Z}}
\def \s{\sigma}
\def \c{\sqrt{{n \choose k}k!}}

\def \x{\sqrt{{n({\delta}) \choose k}k!}}
\def \X{\boldsymbol{X}}
\def \Y{\boldsymbol{Y}}
\def \L{\boldsymbol{L}}
\def \M{\boldsymbol{M}}
\def \N{\boldsymbol{N}}
\def\E{\mathbb{E}}
\def\e{\boldsymbol{E}}
\def\xo{\mathcal{X}}
\def\Eu{\mathcal{E}}
\def\tP{\mathbb{P}}
\def\S{\mathcal{S}}
\def\B{\mathcal{B}}

\renewcommand{\l}[0]{\left }
\renewcommand{\r}[0]{\right}

\def\on{\Omega_{n({\delta})}}
\def\ton{\on^3}
\def\nd{n({\delta})}
\def\ndk{n({\delta_k})}

\def\G{\mathcal{G}}
\def\F{\mathcal{F}}
\def\rhoo{\varrho}
\def\nat{\mathbb{N} \cup \{0\}}
\def\uz{\underline{\z}}
\def\uo{\underline{\o}}
\def\uout{\Upsilon_{\text{out}}}
\def\uoutd{\Upsilon_{\text{out}}^{\D_0}}
\def\el{\mathcal{L}}
\def\els{\el_{\Sigma}}
\def\op{\mathfrak{o}}
\def\out{\mathrm{out}}
\def\inn{\mathrm{in}}
\def\AA{\mathfrak{A}^m}
\def\mm{\mathfrak{m}}
\def\MM{\mathfrak{M}}

\def\wt{\widetilde}
\def\la{\lambda}

\makeatletter
\renewcommand*{\@cite@ofmt}{\hbox}
\makeatother

\begin{document}
\title{Rigidity and Tolerance in Gaussian zeroes and Ginibre eigenvalues: quantitative estimates}
\author{\href{http://math.berkeley.edu/~subhro/}{Subhroshekhar Ghosh}\\ Department of Mathematics\\ University of California, Berkeley\\ subhro@math.berkeley.edu}
\date{}
\maketitle
\begin{abstract}
Let  $\Pi$ be a translation invariant point process on the complex plane $\C$ and let $\D \subset \C$ be a bounded open set whose boundary has zero Lebesgue measure. 
We study the conditional distribution of the points of $\Pi$ inside $\D$ given the points outside $\D$. When $\Pi$ is the Ginibre ensemble or the Gaussian zero process, it been shown in \cite{GP} that this conditional distribution is mutually absolutely continuous with the Lebesgue measure on its support. In this paper, we refine the result in \cite{GP} to show that the conditional density is, roughly speaking, comparable to a squared Vandermonde density. In particular, this shows that even under spatial conditioning, the points exhibit repulsion which is quadratic in their mutual separation.
\end{abstract}

\newpage
\section{Introduction}
\label{intro}
In \cite{GP} the authors studied spatial conditioning in the two main natural examples of repelling point processes on the Euclidean plane, namely the Ginibre ensemble and the Gaussian zero process. They established certain ``rigidity'' phenomena, in the sense that for a bounded open set $\D$ (satisfying some minimal regularity conditions) the points outside $\D$ determine a.s. the number $N$ of points inside $\D$ in the first case and both their number $N$ and their sum $S$ in the second case. Thus, the support of the conditional distribution of the inside points is contained in $\D^N$ in the first case, and in a co-dimension 1 subset of $\D^N$ in the second. It was further established in \cite{GP} that, on this restricted set, the conditional distribution is mutually absolutely continuous with respect to the Lebesgue measure. In this paper, our aim is to refine the result in \cite{GP} and establish quantitative bounds on the conditional density of the inside points (considered as a vector in $\D^N$ by taking them in uniform random order).

We first describe the set up and recall the precise results from \cite{GP}. A point process $\Pi$ on $\C$ is a random locally finite point configuration on the two dimensional Euclidean plane. The probability distribution of the point process $\Pi$ is a probability measure $\P[\Pi]$ on the Polish space of locally finite point sets on the plane. Point processes on the plane have been studied extensively. For a lucid exposition on general point processes one can look at \cite{DV}. The group of translations of $\C$ acts in a natural way on the space of locally finite point configurations on $\C$. Namely, the translation by $z$, denoted by $T_z$, takes the point configuration $\Lambda$ to the configuration $T_z(\Lambda):=\{x+z \big| x \in \Lambda\}$. A point process $\Pi$ is said to be translation invariant if $\Pi$ and $T_z(\Pi)$ have the same distribution for all $z \in \C$. In this work we will focus primarily on translation invariant point processes on $\C$ whose intensities are absolutely continuous with respect to the Lebesgue measure. We will consider simple point processes, namely, those in which no two points are at the same location. A simple point process can also be looked upon as a random discrete measure $[\Pi]=\sum_{z \in \Pi} \del_{z}$.

Let $\D$ be a bounded open set in $\C$ whose boundary has zero Lebesgue measure. Let $\S$ denote the Polish space of locally finite point configurations on $\C$. The decomposition $\D=\D \cup \D^c$ induces a factorization $\S=\S_{\inn} \times \S_{\out}$, where $\S_{\inn}$ and $\S_{\out}$ are respectively the spaces of  finite point configurations on $\D$ and locally finite point configurations on $\D^c$.  This immediately leads to the natural decomposition $\Upsilon=(\Upsilon_{\inn},\Upsilon_{\out})$ for any $\Upsilon \in \S$, and  consequently a decomposition of the point process $\Pi$ as $\Pi=(\Pi_{\inn},\Pi_{\out})$.

For a pair of random variables $(X,Y)$ which has a joint distribution on a product of Polish spaces $\S_1 \times \S_2$, we can define the \textsl{regular conditional distribution} $\gamma$ of $Y$ given $X$ by the family of probability measures  measures $\gamma(s_1,\cdot)$ parametrized by the elements $s_1 \in \S_1$ such that for any Borel sets $A \subset \S_1$ and $B \subset \S_2$ we have 
\[ \P \l( X \in A, Y \in B \r) = \int_{A} \gamma(s_1,B) \, d\,\P[X](s_1) \]
where $\P[X]$ denotes the marginal distribution of $X$.  For details on regular conditional distributions, see, e.g., \cite{Pa}.

Recall that $\S_{\inn}$ and $\S_{\out}$ are Polish spaces. Hence, by abstract nonsense, there exists a \textsl{regular conditional distribution} $\rhoo$ of $\Pi_{\inn}$ given $\Pi_{\out}$. Clearly,  $\rhoo$ can be seen as the distribution of a point process on $\D$ which depends on $\Upsilon_{\out}$.

Let $\uz$ be the vector (of variable length) whose co-ordinates are the points of $\Pi_{\inn}$ taken in uniform random order. We will denote the conditional distribution of $\uz$ given $\Pi_{\out}$ by $\rho$. Formally, it is a family of probability measures $\rho(\uout,\cdot)$ on $\bigcup_{m=0}^{\infty}\D^m$ parametrized by $\uout \in \S_{\out}$.  

For a vector $\underline{\alpha}$, we denote by $\Delta(\underline{\alpha})$ the Vandermonde determinant generated by the co-ordinates of $\alpha$. Note that $|\Delta(\underline{\alpha})|$ is invariant under permutations of the co-ordinates of $\underline{\alpha}$. 

For two measures $\mu_1$ and $\mu_2$ defined on the same measure space $\O$ with $\mu_1 \ll \mu_2$ (meaning $\mu_1$ is absolutely continuous with respect to $\mu_2$), we will denote by $\frac{d \mu_1}{d \mu_2}(\o)$ the Radon Nikodym derivative of $\mu_1$ with respect to $\mu_2$ evaluated at $\o \in \O$.


In Theorems \ref{gingp-1}-\ref{gaf-2} we denote the Ginibre ensemble by $\G$ and the GAF zero ensemble by $\F$. For details on these models, we refer to Section 3 of \cite{GP}, as well as \cite{HKPV}. As before, $\D$ is a bounded open set in $\C$ whose boundary has zero Lebesgue measure. 

In order to provide a concrete framework in which to state the main theorems of this paper (Theorems \ref{gin-2}-\ref{gaf-2}), we quote  Theorems \ref{gingp-1}-\ref{gafgp-2} from \cite{GP}.
\begin{theorem}
 \label{gingp-1}
For the Ginibre ensemble, there is a measurable function $N:\S_{\out} \to \nat$ such that a.s. \[ \text{ Number of points in } \G_{\inn} = N(\G_{\out})\hspace{3 pt}.\]
\end{theorem}

Since a.s. the length of $\uz$ equals $N(\G_{\out})$, we can as well assume that each measure $\rho(\uout,\cdot)$ is supported on $\D^{N(\G_{\out})}$.

\begin{theorem}
\label{gingp-2}
For the Ginibre ensemble, $\P[\G_{\out}]$-a.s. the measure $\rho({\G_{\out}},\cdot)$ and the Lebesgue measure $\el$ on $\D^{N(\G_{\out})}$ are mutually absolutely continuous.
\end{theorem}


\begin{theorem}
 \label{gafgp-1}
For the GAF zero ensemble, \newline
\noindent
(i)There is a measurable function $N:\S_{\out} \to \nat$ such that a.s.  \[ \text{ Number of points in } \F_{\inn} = N(\F_{\out}).\]
(ii)There is a measurable function $S:\S_{\out} \to \C$ such that a.s.  \[ \text{ Sum of the points in } \F_{\inn} = S(\F_{\out}).\]
\end{theorem}

Define the set \[ \Sigma_{S(\F_{\out})} := \{ \uz \in \D^{N(\F_{\out})} : \sum_{j=1}^{N(\F_{\out})} \z_j = S(\F_{\out}) \}\] where $\uz=(\z_1,\cdots,\z_{N(\F_{\out})})$.

Since a.s. the length of $\uz$ equals $N(\F_{\out})$, we can as well assume that each measure $\rho(\uout,\cdot)$ gives us the distribution of a random vector in  $\D^{N({\uout})}$ supported on $\Sigma_{S({\uout})}$.

\begin{theorem}
\label{gafgp-2}
For the GAF zero ensemble, $\P[\F_{\out}]$-a.s. the measure $\rho(\F_{\out},\cdot)$ and the Lebesgue measure $\els$ on $\Sigma_{S(\F_{\out})}$ are mutually absolutely continuous.
\end{theorem}

In this paper, we prove the following quantitative estimates on the density of the conditional distributions with respect to the Lebesgue measure on their support:

\begin{theorem}
 \label{gin-2}
There exist positive quantities $m(\G_{\out})$ and $M(\G_{\out})$ such that a.s. we have
\[m(\G_{\out})|\Delta(\uz)|^2  \le \frac{d \rho(\G_{\out},\cdot)}{d \el}(\uz) \le M(\G_{\out}) |\Delta(\uz)|^2 \]
where $\Delta(\uz)$ is the Vandermonde determinant formed by the co-ordinates of $\uz$ and $\el$ is the Lebesgue measure on $\D^N(\G_{\out})$ . 
\end{theorem}

\begin{theorem}
 \label{gaf-2}
There exist positive quantities $m(\F_{\out})$ and $M(\F_{\out})$ such that a.s. we have
\[m(\F_{\out})|\Delta(\uz)|^2  \le \frac{d \rho(\F_{\out},\cdot)}{d \els}(\uz) \le M(\F_{\out}) |\Delta(\uz)|^2 \]
where $\Delta(\uz)$ is the Vandermonde determinant formed by the co-ordinates of $\uz$ and $\els$ is the Lebesgue measure on $\Sigma_{S(\F_{\out})}$. 
\end{theorem}

Our results show, in particular, that even under spatial conditioning, the points of the Ginibre ensemble or the GAF zero process repel each other at close range, and the quantitative nature of the repulsion (quadratic in the mutual separation of the points) is similar to that of the unconditioned ensembles.

\section{The General Setup}
\label{genset}
We will work in the same abstract setup as in Section 3.3 \cite{GP}. For the reader's convenience, we recount this general setting in this section.

Fix a Euclidean space  $\Eu$ equipped  with a non-negative regular Borel measure $\mu$. Let $\S$ denote the Polish space of  countable locally finite point configurations on $\Eu$.  Endow $\S$ with its canonical topology, namely the topology of convergence on compact sets (which gives $\S$ a canonical Borel $\sigma$-algebra). Fix a bounded open set $\D \subset \Eu$ with $\mu(\partial \D)=0$. Corresponding to the decomposition $\Eu=\D \cup \D^c$, we have ${\S} =\S_{\inn} \times \S_{\out}$, where $\S_{\inn}$ and $\S_{\out}$  denote the spaces of finite point configurations on $\D$ and  locally finite point configurations on $\D^c$ respectively. 

Let $\Xi$ be a measure space equipped with a probability measure $\tP$. For a random variable $Z:\Xi \to \mathcal{X}$ (where $\mathcal{X}$ is a Polish space), we define the push forward $Z_*\tP$ of the measure $\tP$ by  $Z_*\tP(A)=\tP(Z^{-1}(A))$ where $A$ is a Borel set in $\xo$. Also, for a point process $Z':\Xi \to \S$, we can define point processes $Z_{\inn}':\Xi \to \S_{\inn}$ and $Z_{\out}':\Xi \to \S_{\out}$ by restricting the random configuration $Z'$ to $\D$ and $\D^c$ respectively. 

Let $X,X^n: \Xi \rightarrow \S$ be random variables such that $\tP$-a.s., we have $X^n \rightarrow X$ (in the topology of $\S$). We demand that the point processes $X,X^n$ have their first intensity measures absolutely continuous with respect to $\mu$. We can identify $X_{\inn}$ (by taking the points in uniform random order)  with the random vector $\uz$ which lives in $\bigcup_{m=0}^{\infty}\D^m$. The analogous quantity for $X^n$ will be denoted by $\uz^n$.

For our models we can take $\Eu$ to be $\C$, $\mu$ to be the Lebesgue measure, and $\D$ to be a bounded open set whose boundary has zero Lebesgue measure.  

In the case of the Ginibre ensemble, we can define the processes $\G_n$ and $\G$ on the same underlying probability space so that a.s. we have $\G_n \subset \G_{n+1}\subset \G$ for all $n \ge 1$. For reference, see \cite{Go}. We take $(\Xi,\tP)$ to be this underlying probability space, $X^{n}=\G_n$ and $X=\G$.

In the case of the Gaussian zero process, we take $(\Xi, \tP)$ to be a measure space on which we have countably many standard complex Gaussian random variables denoted by $\{ \xi_k \}_{k=0}^{\infty}$. Then $X^n$ is the zero set of the polynomial $f_n(z)=\sum_{k=0}^{n}\xi_k \frac{z^k}{\sqrt{k!}}$, and $X$ is the zero set of the entire function $f(z)=\sum_{k=0}^{\infty}\xi_k \frac{z^k}{\sqrt{k!}}$. The fact that $X^n \to X$ $\tP$-a.s. follows from Rouche's theorem.

\section{Limits of conditional measures: abstract setting}
\label{limcond}
In this section, we introduce an abstract theorem which would enable us to rigorously pass to certain limits of conditional measures. This theorem is a refinement of a corresponding abstract Theorem 6.2 in \cite{GP}, and it enables us deduce the result, namely, the fact that the conditional density is comparable to a squared Vandermonde density.

We recall the following general proposition from \cite{GP}:

\begin{proposition}
\label{meas1}
Let $\Gamma$ be second countable topological space. Let $\Sigma$ be a countable basis of open sets in $\Gamma$ and let $\A:=\{ \cup_{i=1}^k \sigma_i : \sigma_i \in \Sigma, k \ge 1 \}$. Let $c>0$. To verify that two non-negative regular Borel  measures $\mu_1$ and $\mu_2$ on   $\Gamma$ satisfy $\mu_1(B) \le c \mu_2(B)$ for all Borel sets $B$ in $\Gamma$,  it suffices to verify the inequality for all sets in $\A$. 
\end{proposition}

In this Section, we will work in the setup of Section \ref{genset}, specifying $\D$ to be an open ball, and requiring that the first intensity of our point process $X$ is absolutely continuous with respect to the Lebesgue measure on $\Eu$. We further assume that $X$ exhibits rigidity of the number of points. In other words, there is a measurable function $N:\S_{\out} \to \nat$ such that a.s. we have \[\text{ Number of points in } X_{\inn} = N(X_{\out}). \] In such a situation, we can identify $X_{\inn}$ (by taking the points in uniform random order)  with a random vector $\uz$   taking values in $\D^{N(X_{\out})}$.  Studying the conditional distribution $\rhoo(X_{\out},\cdot)$ of $X_{\inn}$ given $X_{\out}$ is then the same as studying the conditional distribution of this random vector given $X_{\out}$. We will denote the latter distribution by $\rho(X_{\out},\cdot)$. Note that it is supported on $\D^{N(X_{\out})}$ (see Section \ref{intro} for details).    

For $m >0$, let $\mathfrak{W}_{\inn}^m$ denote the countable basis for the topology on $\D^m$ formed by open balls contained in $\D^m$ and having rational centres and rational radii. We define the collection of sets $\AA:=\{ \bigcup_{i=1}^k A_i: A_i \in \mathfrak{W}_{\inn}^m , k \ge 1\}$.

Fix an integer $n \ge 0$, a closed annulus $B \subset \D^c$ whose centre is at the origin and which has a rational inradius and a rational outradius, and a collection of $n$ disjoint open balls $B_i$ with rational radii and centres having rational co-ordinates such that $\{ B_i \cap \D^c \}_{i=1}^n \subset B$. Let $\Phi(n,B,B_1, \cdots, B_n)$ be the Borel subset of $\S_{\out}$ defined as follows:
\[ \Phi(n,B,B_1,\cdots,B_n)= \{ \Upsilon \in \S_{\out} : | \Upsilon \cap B |=n, | \Upsilon \cap {\B_i} | =1 \}.\] 
Then the countable collection $\Sigma_{\out}=\{ \Phi(n,B,B_1,\cdots,B_n): n,B, B_i \text{ as above } \}$ is a basis for the topology of $\S_{\out}$. 
Define the collection of sets $\mathcal{B}:=\{\bigcup_{i=1}^k \Phi_i: \Phi_i \in \Sigma_{\out}, k\ge 1 \}$.

We will denote by $\Omega^m$ the event that $|X_{\inn}|=m$, and by $\Omega^m_n$ we will denote the event $\l|X^n_{\inn}\r|=m$. 

\begin{definition}
Let $p$ and $q$ be indices (which take values in   potentially infinite abstract sets), and $\alpha(p,q)$ and $\beta(p,q)$ be non-negative functions of these indices. We write $\alpha(p,q) \asymp_{q} \beta(p,q)$ if there exist positive numbers $k_1(q),k_2(q)$ such that \[ k_1(q)\alpha(p,q) \le \beta(p,q) \le k_2(q)\alpha(p,q) \text{ for all } p,q. \] The main point is that $k_1,k_2$ in the above inequalities are uniform in $p$, that is, all the indices in question other than $q$.  
\end{definition}
We will also use the notation introduced in Section \ref{genset}.
We will define an ``exhausting'' sequence of events as: 
\begin{definition}
A sequence of events $\{\Omega(j)\}_{j \ge 1}$ is said to \textsl{exhaust} another event $\Omega$ if  $\Omega(j) \subset \Omega(j+1) \subset \Omega$  for all $j $ and  $\tP(\Omega \setminus \Omega(j)) \rightarrow 0$ as $j \to \infty$.
\end{definition}

Let $\mathcal{M}(\D^{m})$ denote the space of all probability measures on $\D^m$. For two random variables $U$ and $V$ defined on the same probability space, we say that $U$ is measurable with respect to $V$ if $U$ is measurable with respect to the sigma algebra generated by $V$.  Finally, recall the definition of $\AA$ and $\mathcal{B}$ from the beginning of this section.

Now we are ready to state the following important technical reduction:

\begin{theorem}
\label{abs}
Let $m\ge 0$ be such that $\tP(\Omega^m)>0$. 
Suppose that: 

(a) There is a map $\nu: \Omega^m \rightarrow \mathcal{M}(\D^{m})$ such that for each Borel set $A \subset \D^m$, the random variable $\nu(\cdot,A)$ is measurable with respect to $X_{\out}$. Let $\wt{\nu}:S_{\out}\to \mathcal{M}(\D^m)$ denote the map induced by $\nu$.

(b) For each fixed $j$ we have an event $\Omega(j)$, a sequence $\{ n_k \}_{k \ge 1}$ (which might depend on $j$) and corresponding events $\Omega_{n_k}(j)$    such that: 

(i) Each $\Omega(j)$ is  measurable with respect to $X_{\out}$, and the sequence of events $\Omega(j)$ exhausts $\Omega^m$. 

(ii) $\Omega_{n_k}(j)$ is measurable with respect to $X_{\out}^{n_k}$ and $\Omega_{n_k}(j) \subset \Omega^m_{n_k}$. 

(iii) $\Omega(j) \subset \varliminf_{k \to \infty} \Omega_{n_k}(j) $. 

(iv) For all $A \in \mathcal {A}$ and $B \in \mathcal{B}$ we have
\begin{equation} 
\label{abscond}
\tP [  (\uz^{n_k} \in A) \cap  ( X_{\out}^{n_k} \in B) \cap \Omega_{n_k}(j) ] \asymp _j \int_{(X_{\out}^{n_k})^{-1}(B) \cap \Omega_{n_k}(j)  } \nu({\xi},A) d\tP(\xi) + \vartheta(k;j).   \end{equation}                              
where $\lim_{k \to \infty} \vartheta(k;j)\to 0$ for each fixed $j$.
                                                                                                                                                                         
Then, there are functions $\mm,\MM: \S_{\out} \to \R_+$ such that  a.s. on the event $\Omega^m$ we have 
\begin{equation}
\label{targeteq}
 \mm(\o)\wt{\nu}(\o,A) \le \rho(\o,A) \le  \MM(\o)\wt{\nu}(\o,A)  \end{equation}
for all Borel sets $A$ in $\D^m$, where $\o$ denotes $X_{\out}(\xi), \xi \in \Xi$.
\end{theorem}
\vspace{5 pt}
\begin{remark}
\label{rmabs-1}
The condition that $\Omega(j)$ is  measurable with respect to $X_{\out}$ can be relaxed to the condition that $\widetilde{\Omega}(j)$ is  measurable with respect to $X_{\out}$, for some event $\widetilde{\Omega}(j)$ which satisfies $\P\l( \Omega(j) \Delta \widetilde{\Omega}(j) \r)=0$.
\end{remark}

We defer the proof of Theorem \ref{abs} to Section \ref{techproof}.

We conclude this section with the following simple observation:
\begin{remark}
\label{rmabs-2}
If we have Theorem \ref{abs} for all $m \ge 0$ then we can conclude that (\ref{targeteq}) holds a.e. $\xi \in \Xi$. 
\end{remark}

\section{Reduction from a general $\D$ to a disk}
\label{general}

We will provide a proof that Theorem \ref{gaf-2} in the case where $\D$ is a disk is enough to deduce the case of a general $\D$. The corresponding reduction for Theorem \ref{gin-2} is on similar lines. The first part of this proof is very similar to the corresponding reduction in \cite{GP}. Some additional arguments are introduced subsequently in order to take care of the comparison to the squared Vandermonde.

Let $\D$ be a bounded open set in $\C$ whose boundary has zero Lebesgue measure. By translation invariance of the Ginibre ensemble, we take the origin to be in the interior of $\D$. Let $\D_0$ be a disk (centred at the origin) which contains $\overline {\D}$ in its interior (where $\overline{\D}$ is the closure of $\D$).


Suppose we have the tolerance Theorem \ref{gaf-2} for a disk.  To obtain Theorem \ref{gaf-2} for $\D$ , we appeal to the tolerance Theorem \ref{gaf-2} for the disk $\D_0$. Let the number and the sum of the points in $\D$ be $N$ and $S$ respectively, let $N_0$ and $S_0$ denote the corresponding quantities for the points in $\D_0$.  Define 
\[\Sigma:=\bigg\{(\la_1,\cdots,\la_{N}):\sum_{j=1}^{N}\la_j=S, \la_j \in \D \bigg\} \]
and
\[\Sigma_0 := \bigg\{(\la_1,\cdots,\la_{N_0}):\sum_{j=1}^{N_0}\la_j=S_0, \la_j \in \D_0 \bigg\}.\]


The conditional distribution of the vector of points inside $\D_0$, given $\uoutd$, lives on $\Sigma_0$, in fact it has a density $f_0$ which is positive a.e. with respect to Lebesgue measure on $\Sigma_0$. Let there be $k$ points in $\D_0 \setminus \D$  and let their sum be $s$, clearly we have $N=N_0-k$ and $S=S_0-s$. We parametrize $\Sigma$ by the last $N-1$ co-ordinates. Note that the set $U:=\{(\la_2,\cdots,\la_N): (S-\sum_{j=2}^N \la_j, \la_2,\cdots,\la_N) \in \Sigma\}$ is an open subset of $\D^{N-1}$. Further, we define the set $V:=\{(\la_1,\cdots,\la_k): \la_i \in \D_0 \setminus \D, \sum_{i=1}^k \la_i =s \}$.

Let the points in $\D_0\setminus \D$, taken in uniform random order, form the vector $\underline{\mathbf{z}}=(z_1,\cdots,z_k)$. Then we can condition the vector of points in $\D_0$ to have its last $k$ co-ordinates equal to $\underline{\mathbf{z}}$, to obtain the following formula for the conditional density of the vector of points in $\D$ at $(\z_1,\cdots,\z_N) \in \Sigma$ (with respect to the Lebesgue measure on $\Sigma$):
\begin{equation} 
\label{condng}
f(\z_1,\z_2,\cdots,\z_N)=\frac{f_0(\z_1, \z_2,\cdots,\z_N,z_1,\cdots,z_k)}{\int_U f_0(s-(\sum_{j=2}^N{w_j}), w_2,\cdots,w_N,z_1,\cdots,z_k)dw_2 \cdots dw_N}.
\end{equation}

It is clear that for a.e. $\underline{\mathbf{z}} \in V$, we have $f$ is positive a.e. with respect to Lebesgue measure on $\Sigma$, because the same is true of $f_0$ on $\Sigma_0$.

Let $\uz^0=(\uz,\underline{\mathbf{z}})$.  For estimating $f$ from above and below, recall that a.s. we have \begin{equation} \label{ineq} m(\uoutd)|\Delta(\uz^{0})|^2 \le f_0(\uz^{0}) \le M({\uoutd}) |\Delta(\uz^{0})|^2. \end{equation}

The denominator of the right hand side of (\ref{condng}) does not depend on $\uz$, it only depends on the set $\D$ and $\uout$. Therefore, to obtain the desired upper and lower bounds on $f(\uz)$, it  suffices to estimate $ |\Delta(\uz^{0})|^2 $. 

For two vectors $\underline{\alpha}=(\alpha_1,\cdots,\alpha_m)$ and $\underline{\beta}=(\beta_1,\cdots,\beta_n)$ we define \[\Gamma(\underline{\alpha},\underline{\beta})=\prod_{1\le i \le m, 1\le j \le n} |\alpha_i - \beta_j|.\]  
Then we can write $|\Delta(\uz^{0})|^2=|\Gamma(\uz,\underline{\mathbf{z}})|^2 |\Delta(\underline{\mathbf{z}})|^2 |\Delta(\uz)|^2$. 
Further, since $\partial \D$ is of zero Lebesgue measure, therefore a.s. $\uout \cap \partial \D = \phi = \Upsilon_{\inn} \cap \partial \D$. Hence, there exists $\del(\uout)>0$ such that $\text{ Dist }(\overline{\D},\uout)=\del(\uout)$, in the sense of the distance between a closed set and a compact set. Hence, if $\uz \in \Sigma$,  we have $ \del(\uout) \le |z_i - \z_j| \le 2 \text{ Radius } (\D_0)$ for all $i$ and $j$.

Combining all these observations with equations (\ref{condng}) and (\ref{ineq}) we obtain the desired estimates for $f(\uz)$.

\section{Limiting procedure for the Ginibre ensemble}
\label{ginibre-limit}
As we already saw in Section \ref{general}, it suffices to prove Theorem \ref{gin-2} for $\D$ a disk. To this end, we need to appropriately define the quantities as in the statement of Theorem \ref{abs}, and verify that they satisfy the conditions in that theorem. For the Ginibre ensemble, we define the quantities in exactly the same way as in  Section 8 of \cite{GP}. There we define $\O(j):= \varliminf_{k \to \infty} \O_{n_k}(j)$. In Proposition 8.1 in \cite{GP} it turns out that $\O(j)$ so defined is measurable with respect to $\G_{\out}$ (up to a null set, also see Remark \ref{rmabs-1}) and $\O(j)$-s exhaust $\O^m$, as desired. 

\section{Estimates for finite approximations of the GAF}
\label{estims-gaf}
In this section, we recall several estimates from \cite{GP} which would be necessary to establish Theorem \ref{gaf-2}.  
We focus on the event $\Omega_n^{m,\delta}$ which entails that $f_n$ has exactly $m$ zeroes inside $\D$, and there is a $\delta$ separation between $\partial \D$ and the outside zeroes. The corresponding event for the GAF zero process has positive probability, so by the distributional convergence  $\mathcal{F}_{n} \to \mathcal{F}$, we have that $\Omega_n^{m,\delta}$ has positive probability  (which is bounded away from 0 as $n \to \infty$).

Let us denote the  zeroes of $f_n$ (in uniform random order) inside $\D$  by $\uz=(\zeta_1,\cdots, \zeta_m)$ and those outside $\D$ by $\uo=(\omega_1,\omega_2,\cdots,\omega_{n-m})$. Let $s$ denote the sum of the inside zeroes. Then the conditional density $\rho^n_{\uo}(\uz)$ of $\uz$ given $\uo$ and $s$ is of the following form, supported on the set $\Sigma_s=\{\uz \in \D^m, \sum_{j=1}^m \z_j =s\}$ (see, e.g., \cite{FH}):
\begin{equation}
  \label{condgaf}
\rho^n_{\uo,s}(\uz)=C(\uo,s) \frac{\left| \Delta(\zeta_1,\dots,\zeta_{m},\omega_1,\dots,\omega_{n-m}) \right|^2 }{\displaystyle \left( \sum_{k=0}^n \left| {\sigma_k}{(\uz,\uo)}{\sqrt{{n\choose k} k!} } \right|^2    \right)^{\! \! \! n+1}}
\end{equation}
where $C(\uo)$ is the normalizing factor, and for a vector $\underline{v}=(v_1,\cdots,v_{N})$ we define \[\sigma_k(\underline{v})=\sum_{1 \le i_1 < \cdots <i_k \le N} v_{i_1}\cdots v_{i_k}\] and for two vectors $\underline{u}$ and $\underline{v}$, $\s_k(\underline{u},\underline{v})$ is defined to be $\s_k(\underline{w})$ where the vector $\underline{w}$ is obtained by concatenating the vectors $\underline{u}$ and $\underline{v}$.  

Let $(\uz,\uo)$ and $(\uz',\uo)$ be two vectors of points (under $\mathcal{F}_n$), such that the sum of the co-ordinates of $\uz$ and the same quantity for $\uz'$ are equal to $s$. Then the ratio of the conditional densities at these two vectors is given by
\begin{equation}
 \label{condgaf1}
\frac{\rho^n_{\uo,s}(\uz')}{\rho^n_{\uo,s}(\uz)}= \frac {|\Delta(\uz',\uo)|^2}  {|\Delta(\uz,\uo)|^2}   {\displaystyle \left( \sum_{k=0}^n \left| \frac{\sigma_k(\uz,\uo)}{\sqrt{{n\choose k} k!} } \right|^2    \right)^{\! \! n+1}} \bigg/ {\displaystyle \left( \sum_{k=0}^n \left|\frac {\sigma_k(\uz',\uo)}{\sqrt{{n\choose k} k!} } \right|^2    \right)^{\! \! n+1}}.
\end{equation}
\begin{proposition}
\label{nr2}
On $\Omega_{n}^{m,\delta}$ there are quantities $K(\D,\delta)>0$ and $\X_n(\uo)>0$ such that for any $(\uz,\uz')  \in \D^m \times \D^m$ we have
\[\exp \bigg(- 2mK(\D,\delta)\X_n(\uo) \bigg) \frac{|\Delta(\uz')|^2}{|\Delta(\uz)|^2}\le \frac{|\Delta(\uz',\uo)|^2}{|\Delta(\uz,\uo)|^2} \le \exp \bigg( 2mK(\D,\delta)\X_n(\uo) \bigg) \frac{|\Delta(\uz')|^2}{|\Delta(\uz)|^2}\]
where ${\displaystyle \X_n(\uo) = \l|\sum_{\omega_j \in \mathcal{G}_n \cap \D^c}\frac{1}{\omega_j} \r|+ \l|\sum_{\omega_j \in \mathcal{G}_n \cap \D^c}\frac{1}{\omega_j^2} \r|+ \l( \sum_{\omega_j \in \mathcal{G}_n \cap \D^c}\frac{1}{|\omega_j|^3} \r)}$ and \newline $\E [\X_n(\uo)]\le c_1(\D,m)<\infty .$
\end{proposition}
For the remainder of this section we will restrict $\uz$ and $\uz'$ to lie in the same constant-sum hyperplane. Recall the notation that $s=\sum_{i=1}^{m}\zeta_i$ and \[\Sigma_s := \{\uz \in \D^m : \sum_{i=1}^m \z_i = s \}. \]  Let $D(\uz,\uo)= \left( \sum_{k=0}^n \left| \frac{\sigma_k(\uz,\uo)}{\sqrt{{n\choose k} k!} } \right|^2    \right)$.

When we want to bound the ratio $\l({D(\uz',\uo)} \big/ {D(\uz,\uo)}\r)^{n+1}$ from above and below, it suffices to show that the quantities
\begin{equation}
\label{exp2}
{\displaystyle \l| \sum_{k=0}^n \frac{\overline{\s_k(\uz,\uo)}}{\c} \frac{\s_{k-i}(\uo)}{\c} \r|} \bigg/ {\displaystyle \l( \sum_{k=0}^{n}  \l| \frac{\s_{k}(\uz,\uo)}{\c} \r|^2 \r) }
\end{equation}
and
\begin{equation}
\label{exp3}
{ \displaystyle \l| \sum_{k=0}^n \frac{\overline{\s_{k-i}(\uo)}}{\c} \frac{\s_{k-j}(\uo)}{\c} \r|} \bigg/  {\displaystyle \l( \sum_{k=0}^{n}  \l| \frac{\s_{k}(\uz,\uo)}{\c} \r|^2 \r) }
\end{equation}
for $m \ge i,j \ge 2$ are bounded above by random variables whose typical size is $O(1/n)$ .

The following decomposition of $\s_k(\uz,\uo)$ is simple but useful:
\begin{equation}
\label{exp4}
\s_k(\uz,\uo)=\sum_{r=0}^m \s_r(\uz)\s_{k-r}(\uo).
\end{equation}

The following expansion of $\s_k(\uo)$ in terms of $\s_i(\uz,\uo)$ is more involved:
\begin{proposition}
\label{expansion}
On the event $\O^{m,\del}_n$ we have, for $0 \le k \le n-m $,
\[ \s_k(\uo)=\s_k(\uz,\uo)+\sum_{r=1}^k g_r \s_{k-r}(\uz,\uo)\]
where a.s. the random variables $g_r$ are $O(K(\D,m)^r)$ as $r \to \infty$, for a deterministic quantity $K(\D,m)$ and the constant in  $O$ being deterministic and uniform in $n$ and $\del$.
\end{proposition}
From \cite{GP} equation (22) in Section 9.1.2 we have, for $l \ge m$, 
\begin{equation}
 \label{expansion5}
\frac{\s_{n-l}(\uo)}{\sqrt{{n \choose {n-l}}(n-l)!}}=  \frac{1}{\xi_n} \l[\xi_l + \eta_l^{(n)}   \r] .
\end{equation}
\begin{proposition}
\label{techgaf3}
Let $\eta_l^{(n)}$ be as in (\ref{expansion5}) and $\gamma=\frac{1}{8}$. Then $\exists$ positive random variables $\eta_l$ (independent of $n$) such that a.s. $\l|\eta_l^{(n)}\r|\le \eta_l $, and for fixed $l_0 \in \mathbb{N}$ and large enough $M >0$ we have
 \[(i)\P \l[ \eta_l>\frac{M}{l^{\gamma}} \mathrm{\ for\  some \ } l \ge 1  \r] \le e^{-c_1M^2} \]
 \[(ii)\P \l[ \eta_l>\frac{M}{l^{\gamma}} \mathrm{\ for\  some \ } l \ge l_0  \r] \le e^{-c_2M^2l_0^{\frac{1}{4}}}  \]
where $c_1,c_2$ are constants that depend on the domain $\D$ and on $m$.
\end{proposition}
Define \[ \e_n=\sum_{l=0}^n |\xi_l|^2 , \qquad \L_{ij}^{(n)}=\sum_{l=0}^{n-i \wedge n-j}\frac{\overline{\xi_{l+i}}\xi_{l+j}}{\sqrt{(l+i)_i (l+j)_j}} \hspace{5 pt} , \]
\[
\M_{ij}^{(n)}=\sum_{l=0}^{n-i \wedge n-j}\frac{|\xi_{l+i}|\eta_{l+j}}{\sqrt{ (l+i)_i (l+j)_j}} \hspace{5 pt} , \qquad 
\N_{ij}^{(n)}=\sum_{l=0}^{n-i \wedge n-j}\frac{\eta_{l+i}\eta_{l+j}}{\sqrt{(l+i)_i (l+j)_j}} \hspace{5 pt} .\]
Let \begin{align*} \Y_n= &&  \sum_{i=2}^m \l|\L_{0i}^{(n)} \r| + \sum_{i=2}^m \M_{0i}^{(n)} +  \sum_{i,j \ge 2}^m \l|\L_{ij}^{(n)} \r| + \sum_{i,j \ge 2}^m \M_{ij}^{(n)} + \sum_{i,j \ge 2}^m \M_{ji}^{(n)}   + \sum_{i,j \ge 2}^m \N_{ij}^{(n)}   \hspace{5 pt}. \end{align*}
In \cite{GP} it has been shown (Section 9.1.2 equation (27)) that
\begin{equation}
\label{upperlower}
1-K(m,\D)\frac{\Y_n}{\e_n} \le \frac{D(\uz',\uo)}{D(\uz,\uo)} \le 1+ K(m,\D)\frac{\Y_n}{\e_n} \hspace{5 pt}.
\end{equation}
Regarding $\Y_n$ and $\e_n$, we have the following estimates:
\begin{proposition}
\label{gafestprep}
Given $M>0$ we have:\\
(i) $\P[\Y_n\ge M\log M]\le c(m,\D)/M$ ,\\
(ii) Given $M>0$ there exists $n_0$ such that for $n \ge n_0$ we have $\P[\frac{n}{2}\le|\E_n|\le 2n] \ge 1-\frac{1}{M} \hspace{5 pt}.$ \\
\end{proposition}
These lead to:
\begin{proposition}
\label{gafest}
Given $M>0$ large enough, $\exists n_0$ such that for all $n \ge n_0$ the following is true: with probability $\ge 1- C/M$  we have, on $\Omega_n^{m,\delta}$ , \[ e^{-2K(m,\D)M\log M} \le \l( {D(\uz',\uo)} \bigg/ {D(\uz,\uo)} \r)^{n+1} \le e^{2K(m,\D)M\log M}  \] for all $\uz' \in \Sigma_s$, where $s=\sum_{i=1}^{m}\z_i$ and $(\uz,\uo)$ is randomly generated from $\mathcal{F}_n$.
\end{proposition}
Finally, we arrive at
\begin{proposition}
 \label{unifgaf}
There exist constants $K(m,\D,\delta)$ such that given $M>0$ large enough, we have for $n\ge n_0(m,M,\D)$  the following inequalities hold on $\Omega_n^{m,\delta}$, except for an event of probability $\le c(m,\D)/M$: \[      \exp \bigg({-K(m,\D,\delta)MlogM} \bigg)\frac{|\Delta(\uz'')|^2}{|\Delta(\uz')|^2}  \le \frac{\rho^n_{\uo,s}(\uz'')}{\rho^n_{\uo,s}(\uz')} \le \exp \bigg({K(m,\D,\delta)MlogM}\bigg)\frac{|\Delta(\uz'')|^2}{|\Delta(\uz')|^2}      \]
uniformly for all $(\uz',\uz'') \in \Sigma_s \times \Sigma_s$, where $(s,\uo)$ corresponds to a point configuration picked randomly from $\P\l[\mathcal{F}_n\r]$.
\end{proposition}

The following important estimates deal with inverse power sums of the Gaussian zeroes. 

Let $r_0=$ radius $(\D)$.   Let $\varphi$ be a non-negative radial $C_c^{\infty}$ function supported on $[r_0,3r_0]$ such that $\varphi=1$ on $[r_0+\frac{r_0}{2},2r_0]$ and $\varphi(r_0 + r)=1-\varphi(2r_0+2r)$, for $0\le r \le \frac{1}{2}r_0$. In other words, $\varphi$ is a test function supported on the annulus between $r_0$ and $3r_0$ and its ``ascent'' to 1 is twice as fast as its ``descent''.

\begin{proposition}
\label{invgafest-tail}

(i) The random variables \[ S_l(n):= \int \frac{\widetilde{\varphi}({z})}{z^l} \, d[\F_n](z) + \sum_{j=1}^{\infty} \int \frac{\varphi_{2^j}({z})}{z^l} \, d[\F_n](z) = \sum_{\o \in \F_n \cap \D^c} \frac{1}{\o^l} \quad (\text{ for }l \ge 1) \] and \[
 \widetilde{S}_l(n):= \int \frac{\widetilde{\varphi}({z})}{|z|^l} \, d[\F_n](z) + \sum_{j=1}^{\infty} \int \frac{\varphi_{2^j}({z})}{|z|^l} \, d[\F_n](z) = \sum_{\o \in \F_n \cap \D^c} \frac{1}{|\o|^l} \quad (\text{ for }l \ge 3)\]
have finite first moments which, for every fixed $l$, are bounded above uniformly in $n$. 

(ii) There exists $k_0=k_0(\varphi)\ge 1$,  uniform in $n$ and $l$, such that for $k \ge k_0$ the ``tails'' of $S_l(n)$ and $\widetilde{S}_l(n)$ beyond the disk $2^k \cdot \D $, given by \[\tau_l^{n}(2^k):= \sum_{j=k}^{\infty} \int \frac{\varphi_{2^j}({z})}{z^l} \, d[\F_n](z) \, \,(\text{ for }l \ge 1) \quad \text{ and } \quad \widetilde{\tau}_l^{n}(2^k):= \sum_{j=k}^{\infty} \int \frac{\varphi_{2^j}({z})}{|z|^l} \, d[\F_n](z) \, \,(\text{ for }l \ge 3)\] satisfy the estimates \[ \E \l[ \l| \tau_l^{n}(2^k) \r|\r]\le C_1(\varphi,l)/2^{kl/2} \quad \text{ and } \quad  \E \l[ \l| \widetilde{\tau}_l^{n}(2^k) \r|\r]\le C_2(\varphi,l)/2^{k(l-2)/2}. \]

All of the above remain true when $\F_n$ is replaced by $\F$, for which we use the notations ${S}_l$ , $\widetilde{S}_l$, $\tau_l(2^k)$ and $\wt{\tau}_l(2^k)$ to denote the quantities corresponding to $S_l(n)$ , $\widetilde{S}_l(n)$, $\tau_l^n(2^k)$ and $\tau_l^n(2^k)$.

\end{proposition}

\begin{corollary}
\label{invgafest-tailcor}
For $R=2^k,k\ge k_0$ as in Proposition \ref{invgafest-tail}. We have $\P[|\tau_l^{n}(R)|>R^{-l/4}]\le R^{-l/4}$ and $\P[|\widetilde{\tau}_l^{n}|>R^{-(l-2)/4}]\le R^{-(l-2)/4}$, and these estimates remain true when $f_n$ is replaced with $f$.
\end{corollary}


\begin{proposition}
\label{conv-gaf}
For each $l \ge 1$ we have  $S_l(n) \rightarrow S_l$ in probability, and for each $l \ge 3$ we have $\widetilde{S}_l^n \rightarrow \widetilde{S}_l$ in probability, and hence we have such convergence a.s. along some subsequence, simultaneously for all $l$. 
\end{proposition}

\section{Limiting procedure for GAF zeroes}
\label{limgafz}
In this section, we use the estimates for $\mathcal{F}_n$ to prove Theorem \ref{gaf-2} for a disk $\D$ centred at the origin. We know from Section \ref{general} that this is sufficient in order to obtain Theorems \ref{gaf-1} and \ref{gaf-2} in the general case. We will work in the framework of Section \ref{limcond}. More specifically, we will show that the conditions for Theorem \ref{abs} are satisfied, which will give us the desired conclusion.

\textbf{Notation.} Throughout this section, we will use the notation $\Gamma$ is $O(\delta)$ to imply that the non-negative quantity $\Gamma$ satisfies $\Gamma \le c\delta$ where the positive parameter $\delta \to 0$ and $c$ is a \textbf{universal constant}, which is (crucially) independent of all the other parameters we introduce (like $n$ and $M$). 

In terms of the notation used in Section \ref{limcond}, we have $X^n=\F_n$ and $X=\F$.

\subsection{Overview of the limiting procedure}
\label{ovgafz}
The limiting procedure for the GAF zero process formally involves an appeal to Theorem \ref{abs} - defining the variable $\nu$ and the events $\O(j)$ and $\O_{n_k}(j)$. Further, we need to verify that they indeed satisfy the conditions demanded in Theorem \ref{abs}. The definition of $\nu$ is fairly straightforward, at least in retrospect after the statement of Theorem \ref{gaf-2}.  The main challenge in carrying out the above programme is to simultaneously satisfy two conditions demanded in Theorem \ref{abs}: measurability  of the events with respect to the outside zeroes, and establishing (\ref{abscond}). 

Morally, (\ref{abscond}) corresponds to a statement that for the poynomial ensembles the ratio of the density of the measure $\nu$ with the conditional density of the inside zeroes (given the outside zeroes)  is bounded from above and below in a certain uniform sense. In Section \ref{estims-gaf} we mentioned that in order to obtain such bounds (refer to Proposition \ref{unifgaf}), we need to bound from above expressions such as  $\Y_n$ (refer e.g. to (\ref{upperlower})). However, such expressions involve the coefficients of the Gaussian polynomial, and therefore depend on both the inside and the outside zeroes. Moreover, we would like the criteria defining the events to be somehow consistent, in the sense that having them would imply good bounds (as in (\ref{abscond}) or Proposition \ref{unifgaf}) for all the Gaussian polynomials with large enough degree. This necessitates uniform bounds for the relevant quantities in such a  definition (e.g. the tails of their distributions should be small in some uniform sense), but we only have uniform bounds for the inverse powers of the outside zeroes (refer to Proposition \ref{invgafest-tail}).   

The above factors make the limiting procedure for the Gaussian zero process to be rather involved. We follow the following route: we first define the relevant events (see Definitions \ref{omega} and \ref{omega0} in Section \ref{proofthm}). As is necessary, these definitions make the events measurable with respect to the outside zeroes of $f$ (or $f_n$, as the case may be), but we have to pay the price in that it is not at all clear that these Definitions imply the bounds necessary for obtaining (\ref{abscond}). Indeed, the bounds are not true on the events $\O_{n_k}(j)$ as is, but we will establish that they hold on $\O_{n_k}(j) \setminus E$ where $\P(E)=O(\del_k)$, where $\del_k$ is a decreasing sequence such that $\sum_k \del_k < \infty$. This gives (\ref{abscond}) with the $\vartheta(k;j)$ summand. This deduction can be found in Section \ref{gbounds}. We further need to show that the $\O(j)$-s as defined (see Definition \ref{omega})  exhaust $\O^m$. This is also not clear at all from the definition, and 
has to be deduced in a similar fashion by perturbing the original events, and then taking care of the differences accrued. This is addressed in Section \ref{ELP}. The deductions in Sections \ref{gbounds} and \ref{ELP} will involve several of technical propositions. Their proofs, though straightforward (given the estimates in Section \ref{estims-gaf}), are somewhat peripheral to the broad picture of the limiting procedure. Therefore, these proofs will be deferred to Section \ref{techproofS}.

Our analysis will be driven by two basic parameters: $\del >0$ (to be thought of as small) varying over a decreasing sequence $\del_k$ with $\sum_k \del_k < \infty$, and $M>0$ (to be thought of as large) varying over an increasing sequence $M(j) \uparrow \infty$. We will first fix the parameter $M$ and let $\del \downarrow 0$ along $\del_k$, and then let $M \uparrow \infty$ along $M(j)$. Many of our technical propositions will be true in a regime where $\del$ is small enough depending on $M$, but because of the above order in the taking of limits, this does not create any difficulty.

\subsection{Notations and Choice of Parameters}
\label{parameters}
In this section, we will introduce some notations and define certain parameters, which will be useful for the subsequent analysis. There will be two basic parameters:  $\delta$ and $M$ (see below), and all the other parameters will be defined in terms of these two.  The motivation for these definitions will only be clear when they are invoked at various places in the sections that follow. The reader might omit the detailed definitions on first reading  and refer back when they are used later in the text. The reason we present the definitions together so that the logical dependence between the various parameters become clear, and  there is an organized reference for their future use.    

Let us consider parameters  $\delta>0$ (to be thought of as small), and $M>0$ (to be thought of as large). Let $r_0$ be the radius of the fixed disk $\D$. This $\delta$ bears no relation with the separation parameter $\delta$ considered in Section \ref{estims-gaf}. 

Define, for $j \ge 1$, \[\psi^{f_n}_{2^j,l}=\int \frac{\varphi_{2^j}({z})}{z^l} \, d[\F_n](z)  , \qquad \gamma^{f_n}_{2^j,l}=\int \frac{\varphi_{2^j}({z})}{|z|^l} \, d[\F_n](z),\] \[ \psi^{f}_{2^j,l}=\int \frac{\varphi_{2^j}({z})}{z^l} \, d[\F](z),  \qquad \gamma^{f}_{2^j,l}=\int \frac{\varphi_{2^j}({z})}{|z|^l} \, d[\F](z).\] For $j=0$ replace $\varphi$ by $\wt{\varphi}$ in the above definition.

Define \[ \psi_{l}^{k,f}=\sum_{j=0}^{k}\psi_{2^j,l}^f \hspace{.2 in} \gamma_{l}^{k,f}=\sum_{j=0}^{k}\gamma_{2^j,l}^f \hspace{.2 in} {\psi}_{l}^{f}=\sum_{k=0}^{\infty} \psi_{2^k,l}^{f} \hspace{.2 in} \Psi_{l}^{f}=\sum_{k=0}^{\infty} |\psi_{2^k,l}^{f}|  \hspace{.2 in} {\gamma}_{l}^{f}=\sum_{k=0}^{\infty} \gamma_{2^k,l}^{f} \hspace{3 pt}.\] We define $\begin{displaystyle} \psi_{l}^{k,f_n}, \gamma_{l}^{k,f_n}, {\psi}_{l}^{f_n}, \Psi_{l}^{f_n}  \end{displaystyle}$ and $\gamma_{l}^{f_n}$  as the obvious analogues for $f_n$ (in place of $f$). \newline
Observe that $\begin{displaystyle} {\psi}_{l}^{f_n}=\sum_{\o \in \D^c \cap \F_n} \frac{1}{\o^l} \end{displaystyle}$ and $\begin{displaystyle} {\gamma}_{l}^{f_n}=\sum_{\o \in \D^c \cap \F_n} \frac{1}{|\o|^l} \end{displaystyle}$ .

Let $P_l$ be the $l$-th Newton polynomial expressing the elementary symmetric function of order $l$ in terms of power sums of order $1,2,\cdots,l$. That is, for complex numbers $x_1,\cdots,x_n$, let $s_k=\sum_{j=1}^n x_j^k$ and ${ \displaystyle  e_k=\sum_{i_1 < \cdots < i_k } x_{i_1}x_{i_2}\cdots x_{i_k} }$. Then $e_l=P_l(s_1,\cdots,s_l)$. Note that as a polynomial of $l$ variables, $P_l$  does not depend on $n$ (see \cite{Sta}, Chapter 7).

Let $h:\R^{+} \to \R^{+}$ be a function such that $\sum_{L}^{L+h(L)}k^{-1/8} \to 0$ as $L \to \infty$.  For the sake of definiteness, we set $h(L)=L^{1/16}$. Let $B$ be such that $\sup_{r}\frac{K^r}{(r!)^{1/4}} \le B$, as in Proposition \ref{techgaf3}.

Let $C_{\delta},L_{\delta}$  be  (large) positive integers depending on $\delta$, to be chosen later.

For a non-negative integer $m$ we will consider functions of $C_{\del}$ complex variables $(z_1,\cdots,z_{C_{\del}})$ given by
\[ f_{i,j}^{\del,m} (z_1,\cdots,z_{C_{\delta}})= \l| \sum_{l=1}^{C_{\delta}}\overline{z_{l+i-m}}{z_{l+j-m}} l!  \r|, 2 \le i \le m , 0 \le j \le m \hspace{3 pt}, \]

\[ f_0^{\del,m}(z_1,\cdots,z_{C_{\delta}})=\l(\frac{1}{h(L_{\del})} \sum_{l=L_{\del}}^{L_{\delta}+h(L_{\del})}|z_{l-m}|^2 l!\r) \hspace{3 pt}. \]

Recall the notation that $\Omega^{m,{1}/{M}}$ is the event that there are exactly $m$ zeroes of $f$ in $\D$ and there is a separation of width $\frac{1}{M}$ between the zeroes of $f$ in $\D^c$ and $\partial \D$. Moreover, $\Omega^{m,{1}/{M}}_{n}$ denotes  the analogous event  with $f$ replaced by $f_n$.

With these notations, we make certain choices as follows:

\textbf{I.} We choose an integer $L_{\delta}$ such that \\
(i) $\sum_{k=L_{\delta}}^{L_{\delta}+h(L_{\delta})}k^{-1/8} < \delta^2/C_0$ (with $C_0$ a constant as in the proof of Proposition \ref{compden}).\\
(ii) ${\displaystyle \l| \frac{\sum_{l=L_{\delta}}^{L_{\delta}+h(L_{\delta})}|\xi_l|^2}{h(L_{\delta})}-1 \r| < 1/2 }$ with probability $> 1-\delta^2$.\\
(iii) $\sum_{r=1}^{\infty}\frac{1}{l^{r/16}(r!)^{1/8}} \le \frac{C}{l^{1/16}} <1/2$ for all $l \ge L_{\delta}$, as in the proof of Proposition \ref{compden}.\\
(iv) $L_{\delta}^{1/8 }>K$ where $K=K(\D,m)$ as in Proposition \ref{expansion}.

\textbf{II.} We choose an integer $C_{\delta}$ such that:\\
(i) $C_{\delta}>L_{\delta}+h(L_{\delta}) \hspace{3 pt},$ \\
(ii) $\sum_{l=C_{\delta}}^{\infty}\frac{1}{(l+1)(l+2)} < \delta^6 \hspace{3 pt},$ \\
(iii) $\sum_{l=C_{\delta}}^{\infty} \frac{1}{l^{1/8}\sqrt{(l+1)(l+2)}} <\delta^4 \hspace{3 pt},$\\
(iv) ${\displaystyle e^{-c_2C_{\del}^{1/4}}}<\del^2$, where $c_2$ is as in Proposition \ref{techgaf3}.

\textbf{III.} Let $g(\delta)$ be such that $g(\delta) \to \infty $ as $\delta \to 0$ and \[\P(\Psi_l^{f}>g(\delta)-1)<\delta^2/{C_{\delta}} \hspace{3 pt}, l=1,2,\cdots,C_{\delta} \hspace{3 pt},\] \[\P({\gamma}_3^{f}>g(\delta)-1)<\delta^2 \hspace{3 pt}.\]

\textbf{IV.} Recall that $P_l$ is the $l$-th Newton polynomial expressing the elementary symmetric function of order $l$ in terms of power sums of order $1,2,\cdots,l$. Let $f_{i,j}^{\del,m}$ be the functions defined above. These functions are continuous, and hence uniformly continuous on a compact set.  Let \
$\eps(\delta)$ be such that
\[|f_{i,j}^{\del,m}(P_1(\underline{w_1})),\cdots,P_{C_{\delta}}(\underline{w_1}))-f_{i,j}^{\del,m}(P_1(\underline{w_2})),\cdots,P_{C_{\delta}}(\underline{w_2}))|< \delta^2  \forall ||\underline{w_1}-\underline{w_2}||_{\infty}<\eps({\delta})\] and  $||\underline{w_1}||_{\infty},||\underline{w_2}||_{\infty} \le g(\delta) + 100 $,
for each $2 \le i \le m, 0 \le j \le m$ (here $\underline{w_1},\underline{w_2}$ are  vectors in $\C^{C_{\delta}}$ and $\| \cdot \|$ is the $L^{\infty}$ norm on $\C^{C_{\delta}}$). The same inequality should also hold when $f_{i,j}^{\del,m}$ is replaced by $f_0^{\del,m}$.

\textbf{V.} Let $R_{\delta}>1$ be a large radius (of the form $2^kr_0$ where $r_0$ is radius of $\D$) such that $\sum_{l \ge 1} R_{\del}^{-l/4} < (\eps({\delta})\wedge \delta^2)$. Let $k({\delta})$ be the positive integer such that $R_{\delta}=2^{k({\delta})}r_0$. 

\textbf{VI.} Let $\nd$ be such that \\
(i) Except for an event of probability $< \delta^2$, we have $|\psi_{2^k,l}^{f}-\psi_{2^k,l}^{f_{\nd}}|<\text{min}(\delta^2,\eps(\del))/(k(\delta)+1)$  for all $0 \le k \le k({\delta})$ and all $l \le C_{\delta}$, and $|\gamma_{2^k,3}^{f}-\gamma_{2^k,3}^{f_{\nd}}|<\text{min}(\delta^2,\eps(\del))/(k(\delta)+1)$ for all $0 \le k \le k({\delta})$.\\
(ii) Except for an event of probability $<\delta^2$, we have ${\displaystyle \l| \frac{\sum_{l=1}^{\nd}|\xi_l|^2}{\nd}-1 \r| <1/2}$.\\
(iii) Except for an event of probability $<\delta^2$, we have ${\displaystyle \l| \frac{P_{\inn}^{f_{\nd}}}{|\xi_0|^2}-\frac{P_{\inn}^f}{|\xi_0|^2} \r|<\delta^2}$ where $P_{\inn}^f$ is the product of  zeroes of $f$ inside $\D$ and $P_{\inn}^{f_{\nd}}$ is the analogous quantity for $f_{\nd}$.\\
(iv) Except on an event of probability $ O(\del^2)$, we have

\[ \frac{\displaystyle \l| \sum_{k \le \nd - C_{\del}}\frac{\overline{\sigma_{k-t}(\uo)}\sigma_{k-i}(\uo)}{{\nd \choose k}k!} \r| }{\displaystyle \sum_{k=0}^{\nd}\frac{|\sigma_{k}(\uz,\uo)|^2}{{\nd \choose k}k!}} \le \frac{2\del^2}{\nd}, \hspace{5 pt} t=0,1,2,\cdots,m, \hspace{5 pt} 2\le i\le m \hspace{3 pt}.\]
(v) $\P\l( \Omega^{m,{1}/{M}} \Delta \Omega^{m,{1}/{M}}_{\nd} \r) < \delta^2$.

\begin{remark}
That it is at all possible to find such a $\nd$ is not immediate. Conditions VI.(i), VI.(iii) and VI.(v) are ensured by the convergence of the $\mathcal{F}_n$-s to $\mathcal{F}$ on the compact set ${2^{k({\delta})}}.\D$.  VI.(ii) is covered by the law of large numbers. VI.(iv) is addressed in Proposition \ref{lastestimates1}. Moreover, if we have a sequence $\del_k \downarrow 0$, it is clear that we can choose $\ndk$ such that $\ndk \uparrow \infty$.
\end{remark}

\subsection{Proof of Theorem \ref{gaf-2}}
\label{proofthm}

In this section we provide a proof of Theorem \ref{gaf-2}, based on Theorem \ref{abs}. We will show how to define the relevant events here; the verification of some of the properties demanded in Theorem \ref{abs} requires substantial amount of work, and will be taken up in  subsequent sections.

The notations are from Section \ref{parameters} unless stated otherwise.

\begin{proof} [{{Proof of Theorem \ref{gaf-2}}}]

Following the notation in Theorem \ref{abs}, we begin with  $\O^m$, $m\ge 0$.  The cases $m=0$ and $m=1$ are trivial, so we focus on the case $m \ge 2$. 

Our candidate for $\nu({\xi},\cdot)$ (refer to Theorem \ref{abs}) is the probability measure $Z^{-1}|\Delta(\uz)|^2d\el(\uz)$ on $\Sigma_s$, where $\Delta(\uz)$ is the Vandermonde determinant formed by the coordinates of $\uz $ and $s=S(X_{\out}(\xi))$, $\el$ is the Lebesgue measure on $\Sigma_s$ and $Z$ is the normalizing factor. Here we recall the definition of $S(X_{\out}(\xi))$ from Theorem \ref{gafgp-1} and the definition of $\Sigma_s$ from Section \ref{estims-gaf}. 

We now intend to define the events $\O(j)$ and $\O_{n_k}(j)$ as in Theorem \ref{abs}.

Fix two sequences of positive numbers $M_j \uparrow \infty $ and $\del_k \downarrow 0$ such that $\sum \del_k < \infty$.

Define $z_l=P_l(\psi_{1}^{k({\delta}),f},\cdots,\psi_{l}^{k({\delta}),f})$. 

\begin{definition}
\label{omega}
For any $M, \del >0$ we define the event $\O(M;\del)$ by the following conditions (for clarification of the notations used refer to Section \ref{parameters}):
\begin{enumerate}
\item ${\displaystyle \Omega^{m,{1}/{M}} \text{ occurs} }.$
\item
  \begin{enumerate}
 \item ${\displaystyle |\psi_{s}^{k({\delta}),f}|\le M, \text{ for } s=1,2 }.$
 \item ${\displaystyle |\gamma_{3}^{k({\delta}),f}|\le M }.$
  \end{enumerate}
\item ${\displaystyle f_{i,j}^{\del,m}(z_1,\cdots,z_{C_{\delta}}) \le M, \text{ for } 2 \le i \le m \text{ and } 0 \le j \le m  }.$
\item ${\displaystyle  \frac{1}{M} \le f_0^{\del,m}(z_1,\cdots,z_{C_{\delta}} )\le M }.$
\item ${\displaystyle |\psi_{l}^{k({\delta}),f}|\le g(\del), l=1,\cdots,C_{\del} }.$
\end{enumerate}
Define $\O(j):=\varliminf_{k \to \infty}\Omega(M(j);\del_k)$.
\end{definition}
For any given $\del$, it is clear that $\O(M_j;\del) \subset \O(M_{j+1};\del)$, which implies that $\O(j) \subset \O(j+1)$. The fact that $\O(j)$-s exhaust $\O$ is by no means clear, and will be proved in Section \ref{ELP}. 

Our next goal is to define the sequence of events $\O_{n_k}(j)$ as in Theorem \ref{abs}. 

Corresponding to the sequence $\del_k \downarrow 0$, we have the sequence of positive integers $n_k=n(\del_k) \uparrow \infty$. Given $M,\del>0$, let $z_l':=P_l(\psi_{1}^{k({\delta}),f_{\nd}},\cdots,\psi_{l}^{k({\delta}),f_{\nd}})$.

\begin{definition}
\label{omega0}
We define the event $\O_{n(\del)}^{0}(M)$ as:
\begin{enumerate}
\item ${\displaystyle \Omega^{m,{1}/{M}}_{\nd} \text{ occurs } }.$
\item
   \begin{enumerate}
\item ${\displaystyle |\psi_{s}^{k({\delta}),f_{\nd}}|\le M, \text{ for } s=1,2  }.$
\item ${\displaystyle |\gamma_{3}^{k({\delta}),f_{\nd}}|\le M }.$
   \end{enumerate}
\item ${\displaystyle f_{i,j}^{\del,m}(z_1',\cdots,z_{C_{\delta}}') \le M , \text{ for } 2 \le i \le m, 0 \le j \le m }.$
\item ${\displaystyle \frac{1}{M} \le f_0^{\del,m}(z_1',\cdots,z_{C_{\del}}') \le M  }.$
\item ${\displaystyle |\psi_{l}^{k({\delta}),f_{\nd}}|\le g(\del)+\del^2, \text{ for } l=1,\cdots,C_{\del} }.$
\end{enumerate}
We define the event $\O_{n_k}(j)$ to be $\O_{\ndk}^{0}(M(j)+1)$.
\end{definition}

Clearly, $\O_{n_k}(j)$ is measurable with respect to $X_{\out}^{n_k}$. But it is not immediate that  $\O(j) \subset \varliminf_{k \to \infty} \O_{n_k}(j)$, or that (\ref{abscond}) holds. 
Both of these will be established in Section \ref{gbounds}.

Along with Sections \ref{gbounds} and \ref{ELP}, this shows that Theorem \ref{abs}  establishes Theorem \ref{gaf-2}.
\end{proof}

\subsection{$\O_{n_k}(j)$ is a good sequence of Events}
\label{gbounds}

Our aim in this section is to prove that the events $\O_{n_k}(j)$ defined in Section \ref{proofthm} satisfy the conditions in Theorem \ref{abs}. We state this formally as:
\begin{theorem}
 \label{goodseq}
With definitions as in Section \ref{proofthm}, we have $\O(j) \subset \varliminf_{k \to \infty} \O_{n_k}(j) $, and the $\O_{n_k}(j)$-s satisfy (\ref{abscond}).
\end{theorem}
 
We will establish this through a sequence of propositions. Some of the propositions are technical in nature, and we defer their proofs to Section \ref{techproofS}, so that the main features of the limiting argument are not lost in the details. 

We begin with the event $\Omega(M;\delta)$ and show that, except on an event of small probability, this implies the event $\on^{0} (M+1)$. Heuristically, this means  replacing quantities defined in terms of   $f$ by the corresponding quantities defined in terms of $f_{\nd}$ :

\begin{proposition}
\label{ondef}
For $\del$ small enough (depending on $M$), there exists an event $ E_{\delta}^1$, with $\P(E_{\delta}^1)=O(\delta^2)$, such that $\{ \Omega(M;\delta) \setminus E_{\delta}^1 \} \subset \on^{0} (M+1)$.
\end{proposition}

We defer the proof of Proposition \ref{ondef} to Section \ref{techproof0}.

For the rest of this subsection, the zeroes we discuss are going to be those of $f_{\nd}$, unless otherwise mentioned. $\uz$ and $\uo$ will respectively denote the vectors of the zeroes of $f_{\nd}$ inside and outside $\D$, taken in uniform random order.
In the next proposition prove that in $\sigma_k(\uo)$ is roughly comparable to $\sigma_k(\uz,\uo)$ for a range of $k$ close to $\nd$ but not too close. We set $l=\nd-k$. 

\begin{proposition}
\label{compden}
 Except on an event of probability $O(\del^2)$, we have  $\frac{1}{2}|\sigma_k(\uz,\uo)|\le |\sigma_k(\uo)| \le \frac{3}{2} |\sigma_k(\uz,\uo)|$ for all $ L_{\del}\le l \le L_{\del} + h(L_{\del})$, where $l=\nd -k$ and $h$ is as in Section \ref{parameters}.
\end{proposition}

We defer the proof of Proposition \ref{compden} to Section \ref{techproof1}.

\begin{definition}
 \label{omega1}
We \textbf{ define an event $\on^1(M)$ } by the following conditions, \newline with $z_l''=P_l \l( \psi_1^{f_{\nd}}, \cdots ,\psi_{l}^{f_{\nd}} \r)$ :
\begin{enumerate}
\item ${\displaystyle \Omega^{m,{1}/{M}}_{\nd} \text{ occurs } }.$
\item
   \begin{enumerate}
\item ${\displaystyle |\psi_{s}^{f_{\nd}}|\le M, \text{ for } s=1,2  }.$
\item ${\displaystyle |\gamma_{3}^{f_{\nd}}|\le M }.$
   \end{enumerate}
\item ${\displaystyle f_{i,j}^{\del,m}(z_1'',\cdots,z_{C_{\delta}}'') \le M , \text{ for } 2 \le i \le m \text{ and } 0 \le j \le m }.$
\item ${\displaystyle \frac{1}{M} \le f_0^{\del,m}(z_1'',\cdots,z_{C_{\del}}'') \le M  }.$
\item ${\displaystyle |\psi_{l}^{f_{\nd}}|\le g(\del)+2\del^2, \text{ for } l=1,\cdots,C_{\del} }.$
\end{enumerate}
\end{definition}

It follows from the argument presented after Proposition \ref{compden} that 
\begin{proposition}
\label{compden1}
$ \{ \on^0(M+1) \setminus E_{\delta}^2  \}  \subset  \on^1(M+2)$, for an event $E_{\delta}^2$ with $\P(E_{\delta}^2)=O(\del^2)$ and small enough $\delta$ (depending on $M$). 
\end{proposition}

We defer the proof of Proposition \ref{compden1} to Section \ref{techproof.5}.

The event $\on^1(M)$ enables us to obtain certain  estimates discussed in the following proposition, the benefits of these will be clear subsequently. We mention that $\sigma_k$ here is the same as $\sigma_{k}^{f_{\nd}}$, and $P_{\inn}$ is the same as $P_{\inn}^{f_{\nd}}$.

\begin{proposition}
\label{denest}
Set $z_l''=P_l \l( \psi_1^{f_{\nd}}, \cdots ,\psi_{l}^{f_{\nd}} \r) \hspace{3 pt}.$\\
For all small enough $\del$ (depending on $M$), we have that the following are true on $\on^1(M)$ (except on an event of probability $O(\del^2)$): 
\begin{itemize}
\item (i) \begin{equation} \label{denoequi} \frac{8}{27}\frac{1}{h(L_{\del})}\l(\sum_{L_{\del}}^{L_{\del}+h(L_{\del})}|z_{l-m}''|^2 l! \r) \le \frac{|P_{\inn}^{f_{\nd}}|^2}{|\xi_0|^2}\le 8 \frac{1}{h(L_{\del})} \l(\sum_{L_{\del}}^{L_{\del}+h(L_{\del})} |z_{l-m}''|^2 l! \r). \end{equation}\\
\item (ii) ${\displaystyle \l( \sum_{k=0}^{\nd}\l|\frac{\sigma_k(\uz,\uo)}{\sigma_{\nd}(\uo)} \r|^2(\nd-k)! \r) \ge \frac{1}{10M}\nd} \hspace{3 pt}.$
\end{itemize}
 \end{proposition}

We defer the proof of Proposition \ref{denest} to Section \ref{techproof2}.

Let the bad event in Proposition \ref{denest} (where the conclusions of the proposition do not hold) be denoted by $E_{\delta}^3$.

\begin{definition}
 \label{omega2}
Define $\on^2(M)=\on^1(M)\setminus E_{\delta}^3$.
\end{definition}

\begin{proposition}
\label{lastestimates1}
For all small enough $\delta$, we have, for some constant $c$, \[ \P \l( \frac{\displaystyle  \l| \sum_{k < \nd - C_{\del}}\frac{\overline{\sigma_{k-t}(\uo)}\sigma_{k-i}(\uo)}{{\nd \choose k}k!} \r| }{\displaystyle \sum_{k=0}^{\nd}\frac{|\sigma_{k}(\uz,\uo)|^2}{{\nd \choose k}k!}} \le \frac{2\del^2}{\nd} \r) \ge 1 - c\delta^2 , \hspace{3 pt} t=0,1,2,\cdots,m;\hspace{3 pt} 2\le i\le m \hspace{3 pt}. \]
\end{proposition}

We defer the proof of Proposition \ref{lastestimates1} to Section \ref{techproof3}.

Let the union of all the exceptional events (where the bounds do not hold) in  Proposition \ref{lastestimates1} be denoted by $E_{\delta}^4$ which clearly has probability $O(\delta^2)$. 

\begin{definition}
\label{omega3}
Define $\ton(M)=\on^2(M) \setminus E_{\delta}^4$.
\end{definition}
In the next proposition, we demonstrate that on $\ton(M)$, we indeed have good bounds on ratio of conditional densities along a constant sum submanifold $\Sigma_s$.

\begin{proposition}
 \label{smallchange}
(a) There is an event $E_{\delta}$ with $\P(E_{\delta})=O(\delta)$ such that  $ \{ \Omega(M;\delta) \setminus E_{\delta} \} \subset \ton (M+2) \hspace{3 pt}. $ \\
(b) For a vector of inside and outside zeroes $(\uz,\uo)$ of $f_{\nd}$ such that $\Omega^m_{\nd}$ occurs, let $s$ denote the sum of the co-ordinates of $\uz$ and let $\Sigma_s=\{\uz' \in \D^m : \sum_{j=1}^{m} \z'_j = s\}$. On the event $\ton(M)$, the ratio of conditional densities $\rho^{\nd}_{\uo}(\uz'')/\rho^{\nd}_{\uo}(\uz')$ at any two vectors $\uz'',\uz'$ in $\Sigma_s$ is bounded from above and below by functions of $M$ (uniformly in $\nd$).\\
\end{proposition}

We defer the proof of Proposition \ref{smallchange} to Section \ref{techproof4}.

\begin{remark}
\label{return}
Let $\widetilde{\Omega}(M;\delta)$ denote the event obtained by demanding conditions 1-4 in Definition \ref{omega}. We observe (by appropriately reversing the arguments presented above) that there is an event $E_{\del}'$ with $\P(E_{\del}')=O(\del^2)$ such that \[ \l( \on^3(M+2)  \setminus E_{\del}' \r) \subset \widetilde{\Omega}(M+4;\delta).\] This will be referred to later in Section \ref{ELP}.
\end{remark}

Now we are ready to complete the proof of Theorem \ref{goodseq}
\begin{proof} [\textbf{Proof of Theorem \ref{goodseq}}]

Since $\O(j)= \varliminf_{k \to \infty} \O(M_j;\del_k) $, and  $\O(M_j; \del_k) \subset \O_{\ndk}^0(M_j + 1) = \O_{n_k}(j)$ by Proposition \ref{ondef}, it follows that 
$\O(j) \subset \varliminf_{k \to \infty} \O_{n_k}(j) $. 

To obtain (\ref{abscond}), we introduce some further notations. On the event $\Omega^m_{\nd}$, let $\gamma_{\nd}(s;\uo^{\nd})$ denote the conditional probability measure on the sum $s$ of inside zeroes given the vector of outside zeroes of $f^{\nd}$ to be $\uo^{\nd}$, and $\mu_{\nd}(\uz;s,\uo^{\nd})$ be the conditional measure on the inside zeroes $\uz$ of $f_{\nd}$, given the vector of outside zeroes to be $\uo^{\nd}$ and the sum of the co-ordinates of $\uz$ to be $s$. Let $S$ denote the set of all possible sums of the inside zeroes; clearly $S$ is a bounded open set in $\C$.

From Theorem \ref{largeprob} (which will be proved in Section \ref{ELP}) we have  $\on^0(M+1) \setminus E \subset \ton(M+2)$ where $\P \l(E\r)=O( \delta^2)$, and that for large enough $M$ and small enough $\del$ (depending on $M$), we have $\P \l(\on^0(M+1)\r) \ge \frac{1}{2}\P \l( \Omega^m \r) \ge 100 \del$. 

Hence, heuristically speaking, a `large part' of $\on^0(M+1)$ must be inside $\ton(M+2)$. Put in terms of conditional measures ( conditioning successively on $\uo^{\nd}$ and then on $s$), this leads to the following (for all small enough $\delta$) :

$\exists$ a set $\Omega_{good}$ (measurable with respect to outside zeroes of $f_{\nd}$), $\Omega_{good} \subset \on^0(M+1)$ with 

$\P_{\nd}(\on^0(M+1) \setminus \Omega_{good})< \delta$ such that:

For each $\uo^{\nd} \in \Omega_{good}, \exists$ a measurable set $S_{good}(\uo^{\nd}) \subset S$ with 

$\gamma(S \setminus S_{good}(\uo^{\nd});\uo^{\nd})<\delta$ such that :

For each $s \in S_{good}(\uo^{\nd})$ we have $\mu(\uz;s,\uo^{\nd})(\Sigma_s \setminus H(s)) < \delta$, where $H(s)$ is the set of $\uz \in \Sigma_s$ such that $(\uz,\uo^{\nd}) \in \ton (M+2)$.

We begin by considering \begin{equation} \label{consider} \int_{{B} \cap \Omega_{good}} \l[ \int_{S \cap S_{good}} \l(\int_{A \cap \Sigma_s} \, d\mu_{\nd}(\uz;s,\uo^{\nd}) \r)d\gamma_{\nd}(s;\uo^{\nd}) \r] d\P_{\nd}(\uo^{\nd})  \hspace{3 pt}. \end{equation}
Recall that $\, d\mu_{\nd}$ above  has a density with respect to Lebesgue measure on $\Sigma_s$. Moreover, by definition of $\Omega_{good}$ and $S_{good}$, we have that for $\uo^{\nd}\in \Omega_{good}$ and $s \in S_{good}(\uo^{\nd})$, there exists a $\uz \in \Sigma_s$ satisfying $(\uz,s,\uo^{\nd})\in \ton(M+2)$. As a result, on the submanifold $\Sigma_s$, the conditional probability measure $\, d\mu_{\nd}$ above is proportional to the measure $|\Delta(\uz)|^2 \, d\el(\zeta)$, the constants of proportionality depending on $M$ and being uniform in $\del$ (refer to Proposition \ref{smallchange}) .

Hence we have, \[(\ref{consider})  \asymp_M \int_{{B} \cap \Omega_{good}} \l[ \int_{S \cap S_{good}} \l( \frac{\int_{A\cap \Sigma_s}|\Delta(\uz)|^2\, d\el(\zeta)}{\int_{ \Sigma_s}|\Delta(\uz)|^2\, d\el(\zeta)} \r) d\gamma_{\nd}(s;\uo^{\nd}) \r] d\P_{\nd}(\uo^{\nd}) \hspace{3 pt}.\]

Setting $h(A;s):=\frac{\int_{A\cap \Sigma_s}|\Delta(\uz)|^2\, d\el(\zeta)}{\int_{ \Sigma_s}|\Delta(\uz)|^2\, d\el(\zeta)}$ we have that the last expression is equal to \[ \int_{{B} \cap \Omega_{good}} \int_{S \cap S_{good}} h(A;s) d\gamma_{\nd}(s;\uo^{\nd})d\P_{\nd}(\uo^{\nd}) \hspace{3 pt}. \]

Note that the $h(A;s)$ is bounded between 0 and 1, so at the expense of an additive error of size $O(\delta)$ we have 
\[(\ref{consider}) \asymp_M \int_{{B} \cap \Omega_{good}} \int_{S} h(A;s) d\gamma_{\nd}(s;\uo^{\nd})d\P_{\nd}(\uo^{\nd}) +O(\delta) \]
because $\gamma((S \setminus S_{good});\uo^{\nd})<\delta$ as long as $\uo^{\nd} \in \Omega_{good}$.

The integral above can also be written as \[ \int  h(A,s_{\nd}(\xi)) 1\{ X_{\out}^{\nd}( \xi) \in   B \cap \O_{good}\} (\xi) d{\tP}(\xi) \]
where ${ \displaystyle s_{\nd}(\xi)=\sum_{\z_i \in X_{\inn}^{\nd}(\xi)} \z_i }$ and $1\{E\}(\cdot)$ is the indicator function of an event $E$ \hspace{3 pt}.

Recalling that $\tilde{\P}(\Omega_{\nd}^0(M+1)\setminus \Omega_{good})<\delta$, and doing an argument similar to the previous step, we have, at the expense of another additive error of $O(\delta)$ 
\[(\ref{consider})  \asymp_M \int  h(A,s_{\nd}(\xi)) 1\{ X_{\out}^{\nd}( \xi) \in   B \cap  \O_{\nd}^0(M+1) \} (\xi) d\tilde{\P}(\xi) +O(\delta)\hspace{3 pt}. \]

We could also complete the $S_{good}$ to $S$ and $\Omega_{good}$ to $\Omega^0_{\nd}(M+1)$ in (\ref{consider}) straightaway (at the expense of two additive errors of $O(\delta)$; recall that $\mu_{\nd}$ is a probability measure)   to obtain
the integral (\ref{consider})
\[= \int_{{B} \cap \Omega_{\nd}^0(M+1)} \int_{S} \l(\int_{A \cap \Sigma_s} \, d\mu_{\nd}(\uz;s,\uo^{\nd}) \r)d\gamma_{\nd}(s;\uo^{\nd})d\P_{\nd}(\uo^{\nd}) + O(\delta) \hspace{3 pt}.\]

We state the above relation succinctly as
\begin{equation} \label{final1} {\tP}(A \cap B \cap \Omega^0_{\nd}(M+1))) \asymp_M \int  h(A,s_{\nd}(\xi)) 1\{ X_{\out}^{\nd}( \xi) \in   B \cap  \O_{\nd}^0(M+1) \} (\xi) d{\tP}(\xi) +O(\delta) \hspace{3 pt}.\end{equation}
 
We now set $M=M_j$ and $\del=\del_k$ and note that $\Omega_{\ndk}^0(M_j+1)=\O_{n_k}(j)$.

We note that since we send $k \to \infty$ (and hence $\del_k \to 0$ for fixed $j$), therefore the criteria in the previous Propositions in this section involving $\del$ being small enough depending on $M$ (e.g. as in Propositions \ref{ondef} or \ref{lastestimates1} or \ref{smallchange}) is eventually satisfied. 

The integral on the r.h.s. of (\ref{final1}) is
\[ \int  h(A,s_{n_k}(\xi)) 1\{ X_{\out}^{n_k}( \xi) \in   B \cap \O_{n_k}(j)\} (\xi) d{\tP}(\xi) \hspace{3 pt}.\]  Recall that we can consider the event $A \in \mathcal{A}$ as a subset of $\D^m$, and such sets have piecewise smooth boundary. For the definition of $\mathcal{A}$ refer to Section \ref{limcond}.  For such sets $A$, we have  $h(A,s)$ is continuous in $s$. But if $\o$ denotes $X_{\out}(\xi)$, then $s_{n_k}(\xi) \to S(\o)$ a.s. Hence, $\exists \eps_k \downarrow 0$ such that $\tilde{\P}(|s_{n_k}(\xi)-S(\o)|>\eps_k) \to 0$ as $k \to \infty$. Coupled with the fact that $0 \le h(A,s) \le 1$ this implies that 
\begin{align*}  \int  h(A,s_{n_k}(\xi)) 1\{ X_{\out}^{n_k}( \xi)  \in   B \cap \O_{n_k}(j)\} (\xi) d\tilde{\P}(\xi) && = \int  h(A,S(\o)) 1\{ X_{\out}^{n_k}( \xi) \in   B  \cap \O_{n_k}(j)\} (\xi) d\tilde{\P}(\xi) \\ && + \op_k(1) \end{align*}
where $\op_k(1)$ is a quantity which $\to 0$ uniformly in $j$, as $k \to \infty$.

It remains to observe that \[\int  h(A,S(\o)) 1\{ X_{\out}^{n_k}( \xi) \in   B \cap \Omega_{n_k}(j) \} (\xi) d\tilde{\P}(\xi) = \int_{{B} \cap \Omega_{n_k}(j)}\nu({\xi},A) d\tilde{\P}(\xi)\hspace{3 pt}\] where we recall the definition of $\nu(\xi,\cdot)$ from the proof of Theorem \ref{gaf-2} in Section \ref{proofthm}.

These facts together verify (\ref{abscond}), and completes the proof of Theorem \ref{goodseq}.
\end{proof}

\subsection{$\Omega(j)$-s exhaust $\Omega^m $}
\label{ELP}
In this section we prove that for all small enough $\delta$ (depending on $M$) the event $\O(M;\del)$ (almost) contains an event $\Omega_M$ of  high probability inside $\Omega^m$ (the probability will depend  on $M$) . In this way, we would ensure that $\O(j)=\varliminf_{k \to \infty} \Omega(M_j;\delta_k)$ tends to be of  full measure inside $\Omega^m$ as $M \to \infty$. Recall here that $\Omega^m$ denotes the event that there are exactly $m$ zeroes of $f$ inside $\D$.

In this section the vectors of the inside and  the outside zeroes, denoted respectively by $\uz$ and $\uo$, refer to those of $f_{\nd}$ (unless mentioned otherwise). 

\begin{theorem}
 \label{largeprob}
There exists an event $\Omega_M$ (with $\P(\Omega_M) \to 1$ as $M \to \infty$) such that for each $M$ sufficiently large, we have $ \{ \Omega_M \cap \Omega^{m,{1}/{M}} \setminus F_{\delta} \} \subset \Omega(M+3;\del)$ for all $\del$ small enough (depending on $M$); here $F_{\delta}$ is an event of probability $O(\delta^2)$. In
particular, this implies that ${\displaystyle \O(j)=\varliminf_{\delta_k \to 0} \Omega(M_j;{\delta_k})}$ satisfies $\P \l( \Omega^m \setminus \O(j) \r) \to 0  $ as $j \to \infty$.
\end{theorem}

Before we prove Theorem \ref{largeprob}, we will make some technical observations.

\begin{proposition}
\label{elpt1}
The following criteria imply (for all small enough $\delta$ depending on $M$) that we will be guaranteed \textbf{ condition 3 } defining $\on^0(M+1)$, except on an event $F^1_{\del}$ of probability $O(\del^2)$  :
\begin{itemize}
\item (i) ${\displaystyle \frac{|P_{\inn}^{f_{\nd}}|^2}{|\xi_0|^2}<\sqrt{M}/3 }\hspace{3 pt}.$
\item (ii) ${\displaystyle \l| \sum_{k=\nd-C_\del}^{\nd}\frac{\overline{\sigma_{k-i}(\uo)}\sigma_{k-j}(\uo)}{{\nd \choose k}k!}  \r| \le \frac{\sqrt{M}}{|\xi_{\nd}|^2}; \hspace{3 pt} 0\le i\le m; \hspace{3 pt} 2\le j\le m}\hspace{3 pt}.$
\end{itemize}
\end{proposition}

We defer the proof of Proposition \ref{elpt1} to Section \ref{techproof5}.

\begin{definition}
\label{gdm}
We define an event $\Gamma({\del},M)$ by the following conditions:
\begin{itemize}
\item (a) ${\displaystyle |\psi_{l}^{k({\delta}),f}|\le M-\frac{1}{2}; \hspace{3 pt} l=1,2}  \text{ and }  \gamma_3^{k({\delta}),f} \le M-\frac{1}{2} \hspace{3 pt}.$
\item (b) ${\displaystyle \frac{8}{M}+ \kappa\del^2\le \frac{|P_{\inn}^{f}|^2}{|\xi_0|^2}\le \frac{1}{3}\sqrt{M} - \kappa \del^2  } \text{ where } \kappa >0\text{ is a constant to be declared later }.$
\item (c) ${\displaystyle \l| \sum_{l=0}^{C_{\del}} \frac{\overline{\xi_{l+i}}\xi_{l+j}}{\sqrt{(l+i)_i(l+j)_j}} \r| \le \frac{\sqrt{M }}{4} } ; \hspace{3 pt} 0\le i \le m; \hspace{3 pt} 2 \le j \le m \hspace{3 pt}.$
\item (d) ${\displaystyle \sum_{l=0}^{C_{\del}}\frac{|\xi_{l+i}|}{(l+j)^{1/8}\sqrt{(l+i)_i(l+j)_j}} \le \frac{M^{1/4}}{4} ; \hspace{3 pt} 0 \le i \le m; \hspace{3 pt} 2 \le j \le m }\hspace{3 pt}. $
\item(e) ${\displaystyle |\eta_l^{(\nd)}|\le \frac{M^{1/4}}{l^{1/8}} \text{ for all } l \ge 1 }\hspace{3 pt}.$
\end{itemize}
\end{definition}

With this definition, we have
\begin{proposition}
 \label{elpt2}
$\Gamma({\del},M) \cap \Omega^{m,{1}/{M}} \subset \on^0(M)$ except on an event of probability $O( \del^2)$.
\end{proposition}

We defer the proof of Proposition \ref{elpt2} to Section \ref{techproof6}.

Finally, we define our desired event $\Omega_M$.

\begin{definition}
\label{omegam}
We define the event $\Omega_M$ to be :
\begin{itemize}
\item $(a_0)$ ${\displaystyle |\psi_{l}^{f}|\le M -1 , l=1,2 ; |\gamma_3^{f}| \le M-1} \hspace{3 pt}.$
\item $(b_0)$   ${\displaystyle \frac{8}{M-1}\le \frac{|P_{\inn}^{f}|^2}{|\xi_0|^2}\le \frac{\sqrt{M}}{3}-1 } \hspace{3 pt}.$
\item $(c_0)$ ${\displaystyle \l| \sum_{l=0}^{\infty} \frac{\overline{\xi_{l+i}}\xi_{l+j}}{\sqrt{(l+i)_i(l+j)_j}} \r| \le \frac{\sqrt{M }}{4} -1 } ; \hspace{3 pt} 0\le i \le m; \hspace{3 pt} 2 \le j \le m \hspace{3 pt}. $
\item $(d_0)$ ${\displaystyle \sum_{l=0}^{\infty}\frac{|\xi_{l+i}|}{(l+j)^{1/8}\sqrt{(l+i)_i(l+j)_j}} \le \frac{M^{1/4}}{4}-1 ; \hspace{3 pt} 0 \le i \le m; \hspace{3 pt} 2 \le j \le m } \hspace{3 pt}.$
\item $(e_0)$ ${\displaystyle |\eta_l|< M^{1/4}/{l^{1/8}} \forall l \ge 1 }\hspace{3 pt}.$
\end{itemize}
where $\eta_l$ is as in Proposition \ref{techgaf3}.
\end{definition}

We are now ready to prove Theorem \ref{largeprob}.

\begin{proof} [\textbf{Proof of Theorem \ref{largeprob}}]
 
We claim that for fixed $M$ and small enough $\delta$ (depending on $M$), except on an event of probability $ O( \del^2)$, we have $\Omega_M \subset \Gamma({\del},M)$.
This is going to be true because of the choices of the parameters $C_{\del}$ and $k({\delta})$. 

We begin with the event $\Omega_M$.
First, we look at $(c_0)$ and $(d_0)$. Condition II(ii) in Section \ref{parameters} implies, via a second moment estimate, that the tail $\l| \sum_{l=C_{\delta}}^{\infty} \frac{\overline{\xi_{l+i}}\xi_{l+j}}{\sqrt{(l+i)_i(l+j)_j}} \r|< {\delta}^2$ except with probability $O(\delta^2)$; recall here that we always have $i \vee j \ge 2$. A first moment estimate allows us to draw a similar conclusion about  $\sum_{l=C
_{\del}}^{\infty}\frac{|\xi_{l+i}|}{l^{1/8}\sqrt{(l+i)_i(l+j)_j}}$. Hence, for all small enough $\delta$, we have conditions (c) and (d) defining $\Gamma(\del, M)$ except on an event of probability $O(\del^2)$.

Regarding $(a_0)$, the transition from $\psi_l^{k({\delta}),f}$ to  $\psi_l^f$ takes place via assumption V in Section \ref{parameters} on $R_{\del}$. Recall from Corollary \ref{invgafest-tailcor}) that we have \[ \P(|\tau_{l}(R_{\del})|> R_{\del}^{- l /4}) \le R_{\del}^{- l /4} \hspace{3 pt}.\] Combined with the fact that $R_{\del}^{-1  /4} < \del^2$ (refer to V, Section \ref{parameters}), this enables us to replace $\psi_l^{k({\delta}),f}$ by $\psi_l^{f}$ and $\gamma_3^{k({\delta}),f}$ by $\gamma_3^f$ in $(a_0)$, implying (for all small enough $\del$) condition $(a)$ defining $\Gamma(\del,M)$ except on an event of probability $O(\del)$. For a fixed $M$, we notice that $(b_0)$ implies $(b)$ as soon as $\delta$ is small enough (depending on $M$). Finally, we look at $(e_0)$. Recall from Proposition \ref{techgaf3} that $|\eta_l^{(n)}| \le \eta_l$ a.s. for all $l\ge 1$ in order to deduce $(e)$ in the definition of $\Gamma(\del,M)$ from $(e_0)$ in the definition of $\Omega_M$.

This completes the proof that $\O_M \subset \Gamma (\del,M)$, except for a bad event $F^2_{\del}$ of probability $O(\delta^2)$.

Let $F_{\delta}^3 = F_{\del}^1 \cup F_{\del}^2$, where we recall the definition of $F_{\del}^1$ from the beginning of this proof. Clearly $\P(F_{\delta}^3)=O(\delta^2)$ and we have $\Omega_M \cap \Omega^{m,{1}/{M}} \setminus F_{\delta}^3 \subset \on^0(M+1)$.

However, recall from  Remark \ref{return} that conditions 1-4 defining $\on^0(M+1)$ imply conditions 1-4 defining $\O(M+3;\del)$, except on an event of probability $O(\delta^2)$. Moreover, the complement of condition 5 defining $\O(M+3;\del)$ itself has probability $<\del^2$, by choice of $g({\delta})$ (refer to Section \ref{parameters}). We combine the last two bad events into $F^4_{\delta}$.

Let $F_{\delta}=F^3_{\delta} \cup F^4_{\delta}$. Then we have $\P(F_{\delta})=O(\del^2)$ and $\Omega_M \cap \Omega^{m,{1}/{M}} \setminus F_{\delta} \subset \O(M+3;\del) $, as desired.

Finally, to show that $\P(\Omega_{M}) \to 1$ as $M \to \infty$, we simply observe that $|\psi_{l}^{f}|,\gamma_3^f$ and the random variables appearing in conditions $(c_0)$ and $(d_0)$ are random variables with no mass at $\infty$, and  $\frac{|P_{\inn}^{f}|^2}{|\xi_0|^2}$ does not have an atom at 0; the $\eta_l$-s (appearing in condition $(e_0)$) are taken care of by Proposition \ref{techgaf3} part (i).

We can now take liminf over $\delta_k \to 0$ in $ \{ \Omega_M \cap \Omega^{m,{1}/{M}} \setminus F_{\delta_k} \} \subset \O(M+3;\del_k)$. Recall that  $F_{\delta}$ is of probability $O( \delta^2)$ and $\sum \delta_k < \infty$. Using the Borel Cantelli lemma we thus obtain, for $M$ bigger than some universal constant, \[ \Omega_M \cap \Omega^{m,{1}/{M}} \subset \varliminf_{k \to \infty} \O(M+3;\del_k) \hspace{3 pt}. \] Setting $M=M(j)-3$, this implies  ${\displaystyle \Omega_{M(j)-3} \cap \Omega^{m,{1}/{M(j)-3}}} \subset \Omega(j)$. Since this holds for all large enough $j$, $\P \l( \Omega_M \r) \to 1$  and $\P \l(\Omega^{m} \Delta \Omega^{m,{1}/{M}} \r) \to 0$ as $M \to \infty$, we obtain the fact that $\P(\Omega^m \setminus \O(j)) \to 0$ as $j \to \infty$.
\end{proof}

\subsection{Proofs of some propositions from Sections \ref{gbounds} and \ref{ELP}}
\label{techproofS}
In this section, we include the proofs of the technical propositions which were deferred in Sections \ref{gbounds} and \ref{ELP}.

\subsubsection{Proof of Proposition \ref{ondef}}
\label{techproof0}
This is essentially a consequence of the choices and observations made in Section \ref{parameters}. We also recall  Definitions \ref{omega} and \ref{omega0}.

First notice that by VI.(i), (outside an event of probability $< \delta^2$) for $\Omega(M;\delta)$, we have $|\psi_l^{k({\delta}),f_{\nd}}-\psi_l^{k({\delta}),f}|<\delta^2$ for each $l\le C_{\delta}$, and also $|\gamma_3^{k({\delta}),f_{\nd}}-\gamma_3^{k({\delta}),f}|<\delta^2$; recall here the definitions of   $\psi_l^{k({\delta}),f_{\nd}}$ and $\gamma_3^{k({\delta}),f_{\nd}}$ from Section \ref{parameters}.
Hence we immediately have \textbf{condition 2} in Definition \ref{omega0}; moreover \textbf{ condition 5} in that definition is also clear from condition 5 defining $\O (M;\del)$.

On $\O (M;\del)$, we have $|\psi_{l}^{k({\delta}),f}|<g(\delta)$ for $l\le C_{\del}$. But considering VI.(i), summing over $k \le k({\delta})$ and applying triangle inequality we have $|\psi_{l}^{k({\delta}),f_{\nd}}-\psi_{l}^{k({\delta}),f}|<\eps(\del)$, except for an event of probability $<\del^2$.  Choosing $\underline{w_1}=\{\psi_{l}^{k({\delta}),f}\}_{l=1}^{C_{\del}}$ and $\underline{w_2}=\{\psi_{l}^{k({\delta}),f_{\nd}}\}_{l=1}^{C_{\del}}$ in IV, we have \textbf{condition 3} and \textbf{condition 4} in Definition \ref{omega0} (for all small enough $\delta$, depending on $M$).

\textbf{Condition 1} in Definition \ref{omega0} is the same event as condition 1 in Definition \ref{omega} same except on an event of  probability $<\delta^2$, because of the convergence of $\mathcal{F}_n$ to $\mathcal{F}$ on compact sets, see VI.(v) in Section \ref{parameters}.

Combining all the bad events in the above discussion, we obtain $E_{\delta}^1$, whose probability is $O(\delta^2)$.

\subsubsection{Proof of Proposition \ref{compden}}
\label{techproof1}
We begin by recalling that
\[\frac{\sigma_k(\uo)}{\x}=\frac{\sigma_{k}(\uz,\uo)}{\x} + \sum_{r=1}^{k}g_r \frac{\sigma_{k-r}(\uz,\uo)}{{\nd \choose {k-r}}(k-r)!}\frac{1}{\sqrt{(l+1)\cdots (l+r)}} \]
and ${\displaystyle \frac{\sigma_k(\uz,\uo)}{\x}=\frac{\xi_{l}}{\xi_{\nd}}}$.
As a result, we have
\begin{equation} \label{t} \frac{\sigma_k(\uo)}{\sigma_k(\uz,\uo)}=1+ \sum_{r=1}^{\nd-l}g_r \frac{\xi_{l+r}}{\xi_l}\frac{1}{\sqrt{(l+1)\cdots (l+r)}} \hspace{3 pt}. \end{equation}
The aim is to show that the sum over $r \ge 1$ is small with high probability.
By the A.M.-G.M. inequality we have \[\frac{1}{\sqrt{(l+1)\cdots (l+r)}}\le \frac{1}{l^{r/4}(r!)^{1/4}} \hspace{3 pt},\] and we have already seen in Proposition \ref{expansion} that $|g_r|\le K^r$ where $K$ is a constant that depends on the domain $\D$. As a result, we have
\[\sum_{r=1}^{\nd-l}g_r \frac{\xi_{l+r}}{\xi_l}\frac{1}{\sqrt{(l+1)\cdots (l+r)}} \le \sum_{r=1}^{\nd-l}\frac{K^r}{l^{r/4}(r!)^{1/4}}\l| \frac{\xi_{l+r}}{\xi_l} \r| \hspace{3 pt}.  \]
Since $l\ge L_{\delta}$ is big enough such that $K \le l^{1/8}$ (recall condition I.(iv) from Section \ref{parameters}), we need to estimate $\sum_{r=1}^{\nd-l}\frac{1}{l^{r/8}(r!)^{1/4}}\l| \frac{\xi_{l+r}}{\xi_l} \r|$.
Since $\l| \frac{\xi_{l+r}}{\xi_l} \r|$ satisfies $\P\l(\l| \frac{\xi_{l+r}}{\xi_l} \r| >x \r) \le 1/x^2$, we have $\P\l(\l| \frac{\xi_{l+r}}{\xi_l} \r| >l^{r/16}(r!)^{1/8} \r)\le \frac{1}{l^{r/8}(r!)^{1/4}}$.
If $\l| \frac{\xi_{l+r}}{\xi_l} \r| \le l^{r/16}(r!)^{1/8}$ then for each $k$, the sum over $r \ge 1$ in (\ref{t}) is $\le \sum_{r=1}^{\infty}\frac{1}{l^{r/16}(r!)^{1/8}} \le \frac{C}{l^{1/16}}$ which is $<\frac{1}{2}$ for large enough $l$ (in particular for $l \ge L_{\delta}$), as we desire (recall the definition of $L_{\del}$ from Section \ref{parameters}).
We bound the probability of the complement of this event by a simple union bound over $r$ and see that it is $\le \frac{C_0}{l^{1/8}}$. This gives us a bound for a fixed $l$. Now, we want this to be true with high probability for $L_{\del}\le l \le L_{\del}+h(L_{\del})$. By a union bound of the error probabilities over $l$ in that range, we need to ensure that $\sum_{l=L_{\del}}^{L_{\del}+h(L_{\del})}\frac{1}{l^{1/8}} < \del^2/C_0 $, which we know is true by the choice of $h$ made in Section \ref{parameters}, item I. 

\subsubsection{Proof of Proposition \ref{compden1}}
\label{techproof.5}
Recall Definitions \ref{omega0} and \ref{omega1}. Observe that $\on^1(M)$ differs from $\on^0(M)$ in that $\psi_{l}^{k({\delta}),f_{\nd}}$ is replaced everywhere by $\psi_l^{f_{\nd}}$ (consequently $z_l'=P_l \l( \psi_{1}^{k({\delta}),f_{\nd}},\cdots,\psi_{l}^{k({\delta}),f_{\nd}} \r)$ is replaced by $z_l''=P_l \l( \psi_1^{f_{\nd}}, \cdots ,\psi_{l}^{f_{\nd}} \r)$), $\gamma_3^{k({\delta}),f_{\nd}}$ replaced by $\gamma_3^{f_{\nd}}$ and in condition 5 we have $g(\delta)+\delta^2$ replaced by $g(\delta) + 2\delta^2$. 

On $\on^0(M+1)$, we have $|\psi_{l}^{k({\delta}),f_{\nd}}|\le M+1$ for $l=1,2$ and $|\gamma_{3}^{k({\delta}),f_{\nd}}|\le M+1$. Recall the tail estimates in Corollary \ref{invgafest-tailcor}. Applying these to $R=R_{\del}$ and doing a union bound over  $l$, we get that $|\tau_l(R_{\del})|<\eps(\del) \wedge \delta^2$ for all $l$ and small enough $\del$, and also $|\tilde{\tau}_3(R_{\del})|<\delta^2$, except on an event of probability $O(\del^2)$ (to see this refer to the definition of $R_{\del}$ in Section \ref{parameters} V ). We denote by $(E^2_{\del})'$ the union of the exceptional events in this paragraph, and by $(E^2_{\del})''$ the exceptional event in Proposition \ref{compden}. Define  $E^2_{\del}=(E^2_{\del})' \cup (E^2_{\del})''$. Clearly, $\P \l[E^2_{\delta}\r]=O(\delta^2)$, and $ \{ \on^0(M+1) \setminus E_{\delta}^2  \}  \subset  \on^1(M+2)$. 

\subsubsection{Proof of Proposition \ref{denest}}
\label{techproof2}
(i) We have  ${\displaystyle \frac{\sigma_k(\uz,\uo)}{  \sqrt{ {\nd \choose k}k! } } }=\frac{\xi_{\nd-k}}{\xi_{\nd}}$. Hence, ${\displaystyle \frac{|\sigma_k(\uz,\uo)|^2}{{\nd \choose k} k! |\sigma_{\nd}(\uz,\uo)|^2}}=\frac{1}{\nd!} \frac{|\xi_{\nd-k}|^2}{|\xi_0|^2}.$ \newline
Multiplying both sides by $P_{\inn}$, and observing that ${\displaystyle \frac{P_{\inn}}{\sigma_{\nd}(\uz,\uo)}=\frac{1}{\sigma_{\nd-m}(\uo)}}$ 
we have \begin{equation} \label{qw}  \frac{|\sigma_k(\uz,\uo)|^2}{|\sigma_{\nd-m}(\uo)|^2}(\nd-k)!=|\xi_{\nd-k}|^2\frac{|P_{\inn}|^2}{|\xi_0|^2} \end{equation}
On $\on^1(M)$, we have $\frac{1}{2}|\sigma_k(\uz,\uo)|\le |\sigma_k(\uo)| \le \frac{3}{2} |\sigma_k(\uz,\uo)|$, which means $\frac{2}{3}|\sigma_k(\uo)|\le |\sigma_k(\uz,\uo)| \le 2 |\sigma_k(\uo)|$, for $k=\nd-l$ such that $L_{\del}\le l \le L_{\del}+h(L_{\del})$; refer to Proposition \ref{compden}. Hence
\begin{align*}\frac{4}{9}\frac{1}{h(L_{\del})}\l(\sum_{l=L_{\del}}^{L_{\del}+h(L_{\del})}\frac{|\sigma_k(\uo)|^2}{|\sigma_{\nd-m}(\uo)|^2}(\nd-k)! \r) & \le \frac{1}{h(L_{\del})}\l(\sum_{l=L_{\del}}^{L_{\del}+h(L_{\del})}\frac{|\sigma_k(\uz,\uo)|^2}{|\sigma_{\nd-m}(\uo)|^2}(n-k)! \r) \\ &\le 4 \frac{1}{h(L_{\del})}\l(\sum_{l=L_{\del}}^{L_{\del}+h(L_{\del})}\frac{|\sigma_k(\uo)|^2}{|\sigma_{\nd-m}(\uo)|^2}(\nd-k)! \r)\end{align*}
But by (\ref{qw}), the expression in the centre is ${\displaystyle \frac{|P_{\inn}|^2}{|\xi_0|^2}\frac{1}{h(L_{\del})}\l(\sum_{L_{\del}}^{L_{\del}+h(L_{\del})}|\xi_l|^2 \r)}$.
By our choice of $L_{\del}$, except on an event of probability $O(\del^2)$, we have ${\displaystyle 1/2 < \frac{1}{h(L_{\del})}\l(\sum_{L_{\del}}^{L_{\del}+h(L_{\del})}|\xi_l|^2 \r) < 3/2}$.
Hence
\[  \frac{8}{27}\frac{1}{h(L_{\del})}\l(\sum_{L_{\del}}^{L_{\del}+h(L_{\del})}\frac{|\sigma_k(\uo)|^2}{|\sigma_{\nd-m}(\uo)|^2}(\nd-k)! \r)  \le \frac{|P_{\inn}|^2}{|\xi_0|^2} \]
and
\[ \frac{|P_{\inn}|^2}{|\xi_0|^2}
 \le 8 \frac{1}{h(L_{\del})}\l(\sum_{L_{\del}}^{L_{\del}+h(L_{\del})}\frac{|\sigma_k(\uo)|^2}{|\sigma_{\nd-m}(\uo)|^2}(\nd-k)! \r) \]

which, taken together, is the same as \begin{equation} \label{denoequi} \frac{8}{27}\frac{1}{h(L_{\del})}\l(\sum_{L_{\del}}^{L_{\del}+h(L_{\del})}|z_{l-m}''|^2 l! \r) \le \frac{|P_{\inn}^{f_{\nd}}|^2}{|\xi_0|^2}\le 8 \frac{1}{h(L_{\del})} \l(\sum_{L_{\del}}^{L_{\del}+h(L_{\del})} |z_{l-m}''|^2 l! \r) \end{equation}                                                                                                                                                                                                                                                    which gives us  (i). 

Observe that (\ref{denoequi}) can be re-written as 
\begin{equation}
 \label{denoequi1}
\frac{1}{8} \frac{|P_{\inn}^{f_{\nd}}|^2}{|\xi_0|^2} \le f^{\del,m}_0 \l(z_1'',\cdots,z_{C_{\delta}}'' \r) \le \frac{27}{8} \frac{|P_{\inn}^{f_{\nd}}|^2}{|\xi_0|^2}
\end{equation}

(ii)On $\on^1(M)$, we have $\frac{1}{M } \le f_0^{\del,m}\l(z_1'',\cdots,z_{C_{\delta}}'' \r) \le M$. So, (\ref{denoequi}) implies $ \frac{8}{27M} \le \frac{|P_{\inn}|^2}{|\xi_0|^2}$. Summing (\ref{qw}) from $0$ to $\nd$, noting that by choice of $\nd$ we have $\sum_{l=0}^{\nd} |\xi_l|^2 > \frac{1}{2}\nd$ (except on an event of probability $< \del^2$) and using the above lower bound on ${P_{\inn}}\big/{|\xi_0|^2}$, we get the desired lower bound in part (ii).

\subsubsection{Proof of Proposition \ref{lastestimates1}}
\label{techproof3}
For $\mathcal{F}_{\nd}$, we recall from (\ref{expansion5}) that ${\displaystyle \frac{\sigma_{\nd-l}(\uo)}{\sqrt{{\nd \choose {\nd-l}}(\nd-l)!}}=  \frac{1}{\xi_{\nd}} \l[\xi_l + \eta_l^{(\nd)}   \r] }$, and also recall that ${\displaystyle \sum_{k=0}^{\nd}\frac{|\sigma_{k}(\uz,\uo)|^2}{{\nd \choose k}k!} =\frac{\sum_{l=0}^{\nd}|\xi_l|^2}{|\xi_{\nd}|^2} }$. Substituting these, we get 
\[ \frac{\displaystyle  \l| \sum_{k < \nd - C_{\del}}\frac{\overline{\sigma_{k-i}(\uo)}\sigma_{k-i}(\uo)}{{\nd \choose k}k!} \r| }{\displaystyle \sum_{k=0}^{\nd}\frac{|\sigma_{k}(\uz,\uo)|^2}{{\nd \choose k}k!}} = 
\frac{\displaystyle \l| \sum_{l=C_{\delta}+1}^{\nd} \frac{1}{\sqrt{(l+i)_i(l+j)_j}} \overline{ \l( \xi_{l+i} + \eta_{l+i}^{(\nd)} \r)} \l( \xi_{l+j} + \eta_{l+j}^{(\nd)} \r) \r| }{\displaystyle \sum_{l=0}^{\nd}|\xi_l|^2} \hspace{5 pt}. \]
Set $\e_{\nd}=\sum_{l=0}^{\nd}|\xi_l|^2$.
Expanding the product in each term of the numerator, we observe that it suffices to upper bound the following quantities:
\[{   \l| \sum_{l=C_{\delta}+1}^{\nd} \frac{1}{\sqrt{(l+i)_i(l+j)_j}} \overline{ \xi_{l+i} } \xi_{l+j} \r| } \bigg/ \e_{\nd} \quad \text { and } \quad { \l(   \sum_{l=C_{\delta}+1}^{\nd} \frac{1}{\sqrt{(l+i)_i(l+j)_j}} \l| \xi_{l+i}  \r| | \eta_{l+j}^{\nd} | \r) } \bigg/ \e_{\nd}\]
\[ \text{ and } \l(  \sum_{l=C_{\delta}+1}^{\nd} \frac{1}{\sqrt{(l+i)_i(l+j)_j}} \l| \eta_{l+i}^{\nd}  \r| \l| \eta_{l+j}^{\nd} \r| \r)  \bigg/ \e_{\nd} \hspace{5 pt}.\]
Recall that $(l+i)_i=1$ and $i \vee j \ge 2$.
For the second and the third quantities, we proceed as follows. In Proposition \ref{techgaf3} part (ii) we can set $M=1$ and $l_0=C_{\delta}$ to find that dropping an event of probability $O(\delta^2)$, we have  $|\eta_l^{(\nd)}|<\frac{1}{l^{1/8}}$ for $l \ge C_{\delta}$. On the complement of this small event, we can compute the expectation of the numerator in each case. By our choice of $C_{\delta}$, this will be $  O(\delta^4)$, and by Markov's inequality the numerator is $O(\delta^2)$  with probability $\ge 1-c\delta^2$. For the first quantity to be estimated in the display above, we can derive a similar result by using second moments and the Chebyshev's inequality and exploiting our choice of $C_{\del}$ in Section \ref{parameters}. Moreover, by the law of large numbers, we have that for large enough $\nd$, $\sum_{l=0}^{\nd}|\xi_l|^2 >\frac{1}{2}\nd$ with probability $\ge 1-\delta^2$. Combining all these, we obtain the desired estimate.

\subsubsection{Proof of Proposition \ref{smallchange}}
\label{techproof4}
All zeroes dealt with in this proof are those of $f_{\nd }$.

Part (a) is clear from the arguments prior to Definition \ref{omega3}; we simply take $E_{\del}=\cup_{i=1}^4 E_{\del}^i$\hspace{3 pt}.

Recall that on $\on^3(M)$ we have:\\
\begin{enumerate}
\item
  $
{\displaystyle (a) \hspace{5 pt} {\displaystyle \Omega^{m,{1}/{M}}_{\nd} \text{   occurs} }.} 
{ \hspace {.5 in} (b) \hspace{5 pt}  \displaystyle \gamma_3^{f_{\nd}} \le M .}
{ \hspace {.5 in} (c) \hspace{5 pt}   |\psi_l^{f_{\nd}}|\le M, \text{ for }  l=1,2  .}
  $
\item $f_{i,j}^{\del,m}(z_1'', \cdots , z_{C_{\del}}'') \le M  $ where $z_i''$ are as in Proposition \ref{denest}.

\item ${\displaystyle \l( \sum_{k=0}^{\nd}|\frac{\sigma_k(\uz,\uo)}{\sigma_{n-2}(\uo)}|^2(\nd - k)! \r) \ge \frac{1}{10M}\nd   .}$
\end{enumerate}

Now, recall from Proposition \ref{nr2} and equations (\ref{exp2}) and (\ref{exp3}) that in order to bound the ratio of conditional densities along the submanifold $\Sigma_s$, we need to bound from above:\\
\begin{itemize}
\item \textbf{(a)} The ratio of the Vandermonde determinants $\frac{|\Delta(\uz',\uo)|^2}{|\Delta(\uz,\uo)|^2} \hspace{3 pt}.$
\item \textbf{(b)} $ \frac{\displaystyle  \l| \sum_{k=0}^{\nd }\frac{\overline{\sigma_{k}(\uz,\uo)}\sigma_{k-j}(\uo)}{{\nd \choose k}k!} \r| }{\displaystyle \sum_{k=0}^{\nd}\frac{|\sigma_{k}(\uz,\uo)|^2}{{\nd \choose k}k!}} ; \quad 2\le j\le m \hspace{3 pt}. $
\item \textbf{(c)} $ \frac{\displaystyle  \l| \sum_{k=0}^{\nd }\frac{\overline{\sigma_{k-i}(\uo)}\sigma_{k-j}(\uo)}{{\nd \choose k}k!} \r| }{\displaystyle \sum_{k=0}^{\nd}\frac{|\sigma_{k}(\uz,\uo)|^2}{{\nd \choose k}k!}} ; \quad 0 \le i \le m; \hspace{3 pt} 2\le j\le m \hspace{3 pt}. $
\end{itemize}

As bounds we have:

\textbf{(a)} By definition of $\on^2(M)$, we have $|\psi_1^{f_{\nd}}|,|\psi_2^{f_{\nd}}|$ and $\gamma_3^{f_{\nd}}$ are bounded  by $M$. By Proposition \ref{nr2}, this suffices to upper and lower bound the ratio of Vandermondes in terms of $M$.

\textbf{(b)} In view of (\ref{exp4}), we see that in order to bound (b) from above, it suffices to upper bound (c).

\textbf{(c)} We divide the terms as
\[\frac{\displaystyle \l| \sum_{k=0}^{\nd}\frac{\overline{\sigma_{k-i}}(\uo)\sigma_{k-j}(\uo)}{{\nd \choose k}k!} \r|}{\displaystyle \sum_{k=0}^{\nd}\frac{|\sigma_{k}(\uz,\uo)|^2}{{\nd \choose k}k!}} \le \frac{\displaystyle\l| \sum_{k < n-C_{\delta}}\frac{\overline{\sigma_{k-i}}(\uo)\sigma_{k-j}(\uo)}{{\nd \choose k}k!} \r|}{\displaystyle \sum_{k=0}^{\nd}\frac{|\sigma_{k}(\uz,\uo)|^2}{{\nd \choose k}k!}} + \frac{\displaystyle \l| \sum_{k \ge n-C_{\delta}}\frac{\overline{\sigma_{k-i}}(\uo)\sigma_{k-j}(\uo)}{{\nd \choose k}k!} \r|}{\displaystyle \sum_{k=0}^{\nd}\frac{|\sigma_{k}(\uz,\uo)|^2}{{\nd \choose k}k!}}\]
For the sum over the $k < {\nd}-C_{\delta}$ terms, we appeal to Proposition \ref{lastestimates1} and conclude that the expression corresponding to these terms contributes $\le 2\delta^2/{\nd}$ because $\ton(M) \subset \l( E_{\delta}^4 \r)^c$.

For the expression corresponding to the $k \ge \nd -C_{\del}$ terms, we clear $\nd !$ from numerator and denominator, and divide both by $\l| \s_{\nd -m}(\uo) \r|^2$ (note that $\nd -m$ is the largest $k$ for which $\s_{k}^{\nd}(\uo)$ can be non-zero when $\O_{m,M}^{\nd}$ occurs). We thus reduce ourselves to upper bounding 
\[ \frac{\displaystyle \l| \sum_{k\ge \nd-C_{\del}}\frac{\overline{\sigma_{k-i}(\uo)}\sigma_{k-j}(\uo)}{|\sigma_{\nd-m}(\uo)|^2}(\nd - k)! \r| }{\displaystyle \sum_{k=0}^{\nd} \frac{|\sigma_k(\uz,\uo)|^2}{|\sigma_{\nd-m}(\uo)|^2}(\nd - k)!} \]
From the definitions, it readily follows that $\s_{k-i}(\uo) \big/ \s_{\nd-m}(\uo) = z_q'':= P_{q}(\psi_{1}^{f_{\nd}},\cdots, \psi_{q}^{f_{\nd}}) $ with $q=\nd -k + i -m$.
The numerator above is then simply $f_{i,j}^{\del,m}(z_1'',\cdots,z_{C_{\del}}'')$. We can then combine the bounds in conditions 2 and 3 (laid down at the beginning of this proof) 
to deduce that \[ { \l| \sum_{k\ge \nd-C_{\del}}\frac{\overline{\sigma_{k-i}(\uo)}\sigma_{k-j}(\uo)}{|\sigma_{\nd-m}(\uo)|^2}(\nd - k)! \r| } \bigg/ \l( { \sum_{k=0}^{\nd} \frac{|\sigma_k(\uz,\uo)|^2}{|\sigma_{\nd-m}(\uo)|^2}(\nd - k)!} \r) \le 10M^2/{\nd} \]

In view of (b) and (c) above, we have (also refer to equation (14) in \cite{GP} for greater clarification) \[ 1- 11K(m,\D) M^2/{\nd} \le  \frac{D(\zeta',\o)}{D(\uz,\uo)} \le 1+ 11K(m,\D) M^2/{\nd} \] on $\ton(M)$.    

Recalling (\ref{condgaf1}), we deduce from the above bounds that there exist positive functions $U_1(M)$ and $L_1(M)$ such  that if $(\uz,\uo) \in \ton(M)$ then \[\exp \l( -L_1(M) \r) \le \frac{\rho^{\nd}_{\uo}(\uz')}{\rho^{\nd}_{\uo}(\uz)} \le \exp \l( U_1(M) \r)  \]
where $\uz'$ is another vector in $\Sigma_s$, the constant sum submanifold corresponding to $\uz$. From here, we estimate the ratio at any two points $\uz',\uz''$  belonging to $\Sigma_s$ a.s.:  \begin{equation} \label{referbound}  \exp \l( -L(M) \r) \le \frac{\rho^{\nd}_{\uo}(\uz'')}{\rho^{\nd}_{\uo}(\uz')} \le \exp \l( U(M) \r)  \end{equation}
where $L(M)=U(M)=L_1(M)+U_1(M)$.

\subsubsection{Proof of Proposition \ref{elpt1}}
\label{techproof5}
First recall that by definition of $\nd$, we have ${\displaystyle \frac{1}{2}\le \frac{\sum_{l=0}^{\nd}|\xi_l|^2 }{\nd}\le \frac{3}{2}}$ except on a set of probability $<\del^2$. Also recall the notation 
\[z_l''=P_l(\psi_1^{f_{\nd}},\cdots,\psi_{C_{\delta}}^{f_{\nd}})=\sum_{\o_{i_j} \ne \o_{i_k} \ \in \D^c \cap \F_{\nd}}\frac{1}{\o_{i_1} \cdots \o_{i_l}} \hspace{3 pt}.\] 

Moreover, we have ${\displaystyle \frac{|\sigma_k(\uz,\uo)|^2}{{\nd \choose k}k!}=\bigg| \frac{\xi_{\nd-k}}{\xi_{\nd}}} \bigg|^2$. As a result, condition (ii) implies
\[ \frac{ \l| \sum_{k=\nd-C_\del}^{\nd}\frac{\overline{\sigma_{k-i}(\uo)}\sigma_{k-j}(\uo)}{{\nd \choose k}k!}  \r|}{\sum_{k=0}^{\nd}\frac{|\sigma_k(\uz,\uo)|^2}{{\nd \choose k}k!}}\le 2\sqrt{M}/{\nd} \hspace{3 pt}. \]
On the left hand side, we can replace $\nd!$ by $\l|\sigma_{\nd-m}(\uo)\r|^2$ both in the numerator and the denominator, and combined with the observation that $\sigma_k(\uo)/{\sigma_{\nd-m}(\uo)} = z_{\nd-k-m}''$ this would lead to
\[\frac{\l| \sum_{l=0}^{C_{\del}}\overline{z_{l+i-m}''} z_{l+j-m}''l! \r|}{\sum_{k=0}^{\nd}\frac{|\sigma_k(\uz,\uo)|^2}{|\sigma_{\nd-m}(\uo)|^2}(\nd-k)!}\le 2\sqrt{M}/{\nd}\]
where we switched from $k$ to the new variable $l=\nd - k$. Observe that the numerator is simply $f_{i,j}^{\del,m}(z_1'', \cdots, z_{C_{\del}}'')$
However, summing (\ref{qw}) from $0$ to $\nd$ and noting that by choice of $\nd$ we have $\sum_{l=0}^{\nd} |\xi_l|^2 < \frac{3}{2}{\nd}$ (except on  an event of probability $< \del^2$), we have
\[\sum_{k=0}^{\nd}\frac{|\sigma_k(\uz,\uo)|^2}{|\sigma_{\nd-m}(\uo)|^2}(\nd-k)! \le \frac{3}{2}\nd \frac{|P_{\inn}|^2}{|\xi_0|^2} \hspace{3 pt}.\]
Using this and condition (i) above, we have $f_{i,j}^{\del,m} \le M$. The choice of $k(\delta)$ and the tail estimates in Proposition \ref{invgafest-tail} and Corollary \ref{invgafest-tailcor}  enable us to replace $\psi_l^{f_{\nd}}$ by $\psi_l^{k({\delta}),f_{\nd}}$, and $z_l''$ by $z_l'=P_l(\psi_1^{k({\delta}),f_{\nd}},\cdots,\psi_{C_{\delta}}^{k({\delta}),f_{\nd}})$, with an additive error of size $O(\delta^2)$, except on an event of probability $O(\delta^2)$. This is on similar lines to the switch from $z_l'$ to $z_l''$ we did in defining $\on^1(M)$. Tallying all these, we obtain \textbf{ condition 3 }  on $\on^0(M+1)$, for all small enough $\delta$ depending on $M$.

\subsubsection{Proof of Proposition \ref{elpt2}}
\label{techproof6}
First recall Definitions \ref{gdm} and \ref{omega0}. We start with $\Gamma({\del},M) \cap \Omega^{m,{1}/{M}} $. The choice of $\nd$ (condition VI.(v)) in Section \ref{parameters} ensures  that except on an event of probability $< \delta^2$ , we have that $\Omega^{m,{1}/{M}}_{\nd}$ occurs, thereby verifying \textbf{condition 1} in Definition \ref{omega0}.  By choice of $k({\delta})$ and $\nd$, criterion  $(a)$ defining $\Gamma(\del,M)$ ensures \textbf{condition 2} in Definition \ref{omega0}, except on an event of probability $O(\del^2)$ for $\del$ small enough (depending on $M$). The complement of \textbf{condition 5} in Definition \ref{omega0} is itself an event of probability $O( \del^2)$, by choice of $g(\del),k({\delta})$ and $n({\del})$.

The choice of $\nd$ ensures that $ \l| {|P_{\inn}^{f}|^2} \big/ {|\xi_0|^2} - {|P_{\inn}^{f_{\nd}}|^2} \big/ {|\xi_0|^2} \r|< \del^2$ except on an event of probability $O(\del^2)$. Recall condition (b) defining $\Gamma(\delta,M)$  and equation (\ref{denoequi1}). Let us recall the notations $z_l''=P_l(\psi_1^{f_{\nd}},\cdots, \psi_l^{f_{\nd}})$ and $z_l'=P_l(\psi_1^{k(\del),f_{\nd}},\cdots, \psi_l^{k(\del),f_{\nd}})$ The choice of  $k({\delta})$ (condition V) in Section \ref{parameters} enables us to replace $f_0^{\del,m}(z_1'',\cdots,z_{C_{\del}}'')$ by $f_0^{\del,m}(z_1',\cdots,z_{C_{\del}}')$, on similar lines to the switch from $z_l'$ to $z_l''$ in the definition of $\on^1(M)$. The upshot of all this is that except on an event of probability $O(\delta^2)$ we have that on $\Gamma(\del,M)$ the following is true:
\[ f_0^{\del,m}(z_1',\cdots,z_{C_{\del}}') \ge f_0^{\del,m}(z_1'',\cdots,z_{C_{\del}}'') - \l| f_0^{\del,m}(z_1',\cdots,z_{C_{\del}}') - f_0^{\del,m}(z_1'',\cdots,z_{C_{\del}}'') \r| \]
\[ \ge \frac{1}{8} \frac{|P_{\inn}^{f_{\nd}}|^2}{|\xi_0|^2} - \l| f_0^{\del,m}(z_1',\cdots,z_{C_{\del}}') - f_0^{\del,m}(z_1'',\cdots,z_{C_{\del}}'') \r| \]
\[ \ge \frac{1}{8} \frac{|P_{\inn}^{f}|^2}{|\xi_0|^2} - \frac{1}{8} \l| \frac{|P_{\inn}^{f_{\nd}}|^2}{|\xi_0|^2} - \frac{|P_{\inn}^{f}|^2}{|\xi_0|^2} \r| - \l| f_0^{\del,m}(z_1',\cdots,z_{C_{\del}}') - f_0^{\del,m}(z_1'',\cdots,z_{C_{\del}}'') \r| \]
\begin{equation} \label{exh1} \ge \frac{1}{8}\frac{|P_{\inn}^{f}|^2}{|\xi_0|^2} - \kappa_1 \del^2 \end{equation}
Similarly except on an event of probability $O(\delta^2)$,  we have that on $\Gamma(\del,M)$ the following is true:
\begin{equation} \label{exh2} f_0^{\del,m}(z_1',\cdots,z_{C_{\del}}') \le \frac{27}{8}\frac{|P_{\inn}^{f}|^2}{|\xi_0|^2} + \kappa_2 \delta^2  \end{equation}
Here $\kappa_1$ and $\kappa_2$ are two positive constants. In condition (b) defining $\Gamma(\del,M)$, we choose $\kappa$ to be $\text{ max }(\kappa_1,\kappa_2,1)$. The last pair of inequalities (\ref{exh1}) and (\ref{exh2}) gives us \textbf{condition 4} in Definition \ref{omega0}, for all $M$ sufficiently large (in fact, as soon as $\sqrt{M} \ge 9/8$).

Finally, we have already seen how criteria (i) and (ii) in Proposition \ref{elpt1} suffice to imply \textbf{condition 3} in Definition \ref{omega0}, except on an event of probability $O(\del^2)$. Criterion (i) is implied by the upper bound in condition $(b)$ defining $\Gamma(\delta,M)$, applied along with the inequality $\l| \frac{|P_{\inn}^{f_{\nd}}|^2}{|\xi_0|^2} - \frac{|P_{\inn}^{f}|^2}{|\xi_0|^2} \r|<\del^2$ (which holds everywhere except on an event of probability $<\del^2$ by definition of $\nd$).

Thus, it remains to show that conditions labeled by $(c)$,$(d)$ and $(e)$ above together imply criterion (ii) in Proposition \ref{elpt1}.
To this end, we refer to Section \ref{estims-gaf} where we obtain expressions of $\sigma_k(\uo)$ as expansions in ${\xi_l} \big/ {\xi_{\nd}}$ (as in (\ref{expansion5})). 
From such expansions, we find that in order to establish criterion (ii), it suffices to show
\[\frac{1}{|\xi_{\nd}|^2}\l| \sum_{l=0}^{C_{\del}} \l(\frac{\overline{\xi}_{l+i} + \overline{\eta}_{l+i}^{(\nd)} }{\sqrt{(l+i)_i}}  \r)  \l(\frac{\xi_{l+j}+ \eta_{l+j}^{(\nd)}}{\sqrt{(l+j)_j}}    \r)   \r| \le \frac{\sqrt{M}}{|\xi_{\nd}|^2} \hspace{3 pt}.\]
Recall from Proposition \ref{techgaf3} that except with probability $O(\del^2)$, we have for all $l\ge 1, |\eta_l^{(\nd)} |\le |\eta_l| \le l^{-1/8}$.
It is then a straightforward calculation using the triangle inequality that indeed (c), (d) and (e) defining $\Gamma(\del,M)$ imply the above inequality, for $0\le i \le m; 2 \le j \le m$. 

\section{Proof of Theorem \ref{abs}}
\label{techproof}
We claim that it suffices to show that for every $j_0$ , we have for any $A \in \mathcal{A}$ and any Borel set $B$ in $\S_{\out}$
\begin{equation}
\label{abs1} 
      \tP \l( (\uz(X_{\inn}) \in A) \cap (X_{\out} \in B)  \cap \Omega(j_0) \r) \asymp_{j_0} \int_{X_{\out}^{-1}(B) \cap \Omega(j_0) } \nu({\xi},A) d{\tP}(\xi) \hspace{3 pt}.
\end{equation}
First, note that (\ref{abs1})  implies the conclusion (\ref{targeteq}) of Theorem \ref{abs}   a.s.  on the event $\Omega(j_0)$ with the quantities $m$ and $M$ depending only on $j_0$. To see this, let $\xo$ be the sigma-algebra generated by $X_{\out}$ , and let $\xo(j_0)$ be the collection of sets formed by intersecting the sets of $\xo$ with $\Omega(j_0)$. So, every set in $\xo(j_0)$ is of the form $X_{\out}^{-1}(B)\cap \Omega(j_0)$ for some Borel set $B \subset \S_{\out}$. We can then consider the finite measure space $(\Omega(j_0),\xo(j_0),\tP)$.  For fixed $ A \in \mathcal{A}$ and any Borel set $ B $  in $\S_{\out}$  we have 
\begin{equation} \label{abs2}
 \tP \l( (\uz(X_{\inn}) \in A) \cap (X_{\out} \in B) \cap \Omega(j_0)  \r) = \int_{X_{\out}^{-1}(B) \cap \Omega(j_0) } \rho(X_{\out}(\xi),A) d\tP(\xi)  \hspace{3 pt}.
\end{equation}
We can now compare the (\ref{abs1}) and (\ref{abs2}) for fixed $A \in \mathcal{A}$ and all  Borel sets $B$  in $\S_{\out}$. Considering this on the finite measure space $(\Omega(j_0),\xo(j_0),\tP)$, we obtain   obtain (\ref{targeteq}) for  a.s. $\xi \in \Omega(j_0)$ and a fixed $A \in \mathcal{A}$. The constants $m$ and $M$ obtained are the same as those appearing in (\ref{abs1}) for this $j_0$.  Here we use the fact that all the sides in (\ref{targeteq}) are measurable with respect to $X_{\out}$.
We conclude that  for  a.s. $\xi \in \Omega(j_0)$, (\ref{targeteq}) is simultaneously true for all sets  $A \in \mathcal{A}$, because $\mathcal{A}$ is a countable collection of sets.
From here, we  use Proposition \ref{meas1} to obtain (\ref{targeteq}) for all Borel sets $A$ in $\D^m$ and  a.s. $\xi \in \Omega(j_0)$.
But the $\Omega(j)$-s exhaust $\Omega^m$, so we have the conclusion of Theorem \ref{abs} holding a.e.   $ \xi \in \Omega^m$. The constants however now depend on $\o=X_{\out}(\xi)$, because  for any given $\xi$, we simply take the constants coming out of (\ref{abs1}) for the minimal $j_0$ for which $\xi \in \Omega(j_0)$; recall here that $\Omega(j_0)$ is measurable with respect to $X_{\out}$ .

Now we focus on proving (\ref{abs1}). For any Borel set $B$ in $\S_{\out}$,  given $\eps>0$ we can find a  $\b \in \mathcal{B}$  such that  $\tP \l( \l( X_{\out}^{-1}(B) \cap \Omega(j_0) \r) \Delta X_{\out}^{-1} (\b) \r)< \eps$. This can be seen by considering the push forward probability measure $(X_{\out})_*\tP$ on $\S_{\out}$; recall that $\Omega(j_0)$ is also measurable with respect to $X_{\out}$.
The aim of this reduction is to exploit the fact that as $k\to \infty$, we have $1_{\b}\l({X_{\out}^{n_k}}\r) \rightarrow 1_{\b}\l({X_{\out}}\r)$ a.s. 

We start from (\ref{abscond}) applied to $A,\b,j_0$ (where $A \in \mathcal{A}$ ) and intend to derive (\ref{abs1}) for $A,B,j_0$. We want to show that both sides   of (\ref{abscond}) converge to the appropriate quantities in (\ref{abs1}) as $k \to \infty$ and $ \eps \to 0$. 

In what follows, we will denote by $\tP_h$ the non-negative finite measure on $\Xi$ obtained by setting $d\tP_h(\xi)=h(\xi)d\tP(\xi)$ where $h:\Xi \to [0,1]$ is a measurable function. We note that for any event $E$, we have $0 \le \tP_h(E) \le \tP(E) \le 1 $. For the rest of this proof, the value of $j_0$ is  held fixed.

We start with $\tP_h [  (X_{\inn}^{n_k} \in A) \cap  ( X_{\out}^{n_k} \in \b) \cap \Omega_{n_k}(j_0) ]$. This is equal to \[ \tP_h \bigg[  (X_{\inn} \in A) \cap  ( X_{\out} \in \b) \cap \Omega_{n_k}(j_0) \bigg] + o_k(1;\b) \] where $o_k(1;\b)$ stands for a quantity that tends to $0$ as $k \to \infty$ for fixed $\b$. This step uses the fact that $1_{\b}\l(X^{n_k}_{\out}(\xi)\r) \to 1_{\b}\l(X_{\out}(\xi)\r)$ a.s. The last expression above equals
\[  \tP_h \bigg[  (X_{\inn} \in A) \cap  \l( ( X_{\out} \in B) \cap \Omega(j_0) \r) \cap \Omega_{n_k}(j_0) \bigg]  + o_{\eps}(1) + o_k(1;\b) \]  where $o_{\eps}(1)$ denotes a quantity that tends to $0$ uniformly in $k$ as $\eps \to 0$.  Observe that  $(X_{\out})^{-1}(B) \cap \Omega(j_0) \subset \Omega({j_0}) \subset \varliminf_{k \to \infty} \Omega_{n_k}(j_0)$ and $\tP \l( \varliminf_{k \to \infty} \Omega_{n_k}(j_0) \Delta \l( \bigcap_{l \ge k} \Omega_{n_l}(j_0)  \r)  \r) = \op_k(1)$ where $\op_k(1)$ denotes a quantity such that $\op_k(1) \to 0$ as $k \to \infty$, uniformly in $B$ and $\eps$. Hence we have \[ \tP_h \bigg[(X_{\inn}^{} \in A) \cap ( X_{\out} \in B) \cap \Omega(j_0) \cap \Omega_{n_k}(j_0)   \bigg] \]\[ = \tP_h \bigg[  (X_{\inn} \in A) \cap  ( X_{\out} \in B) \cap \Omega(j_0) \cap \l( \bigcap_{l \ge k} \Omega_{n_l}(j_0)  \r) \cap \Omega_{n_k}(j_0) \bigg]  + \op_k(1)  \hspace{3 pt}.  \]
But $\l( \bigcap_{l \ge k} \Omega_{n_l}(j_0) \r) \cap \Omega_{n_k}(j_0) = \l( \bigcap_{l \ge k} \Omega_{n_l}(j_0) \r) $ and 
\[ \tP_h \bigg[  (X_{\inn} \in A) \cap  ( X_{\out} \in B) \cap \Omega(j_0) \cap \l( \bigcap_{l \ge k} \Omega_{n_l}(j_0) \r) \bigg] = \tP_h \bigg[  (X_{\inn} \in A) \cap  ( X_{\out} \in B) \cap \Omega(j_0) \cap \varliminf_k \Omega_{n_k}(j_0)  \bigg] \]\[ + \op_k(1).  \]
But again,  $( X_{\out} \in B ) \cap \Omega(j_0) \subset \varliminf_k \Omega_{n_k}(j_0)$, so the upshot of all this is that \[\tP_h \bigg[  (X_{\inn}^{n_k} \in A) \cap  ( X_{\out}^{n_k} \in \b) \cap \Omega_{n_k}(j_0) \bigg] = \tP_h \bigg[  (X_{\inn} \in A) \cap  ( X_{\out} \in B) \cap \Omega(j_0) \bigg] + o_k(1;\b) + \op_k(1) + o_{\eps}(1) . \]
We apply the above reduction  to the left hand side of (\ref{abscond}) with $h \equiv 1$ and to the right hand side of (\ref{abscond}) with $h(\xi)=\nu(\xi,A)$, and obtain

\[\tP \bigg[  (X_{\inn} \in A) \cap  ( X_{\out} \in B) \cap \Omega(j_0) \bigg]  \asymp_{j_0}  \l( \int_{(X_{\out})^{-1}(B) \cap \Omega(j_0)} \nu({\xi},A) d\tP(\xi) \r) + o_{\eps}(1) + o_k(1;\b) + \op_k(1)\]\[ + \vartheta(k;j_0) \hspace{3 pt}. \]

Recall from the statement of Theorem \ref{abs}  that $\vartheta(k;j_0) \to 0$ as $k \to \infty$.
First letting $k \to \infty$ in the above (with $\b$ held fixed) and then $\eps \to 0$ (i.e., letting $\b \to B$ in the sense that $\P(\b \Delta B) \to 0$) in the above, we obtain (\ref{abs1}).  

This completes the proof of Theorem \ref{abs}.

\textbf{Acknowledgements.} We are very grateful to Fedor Nazarov for suggesting the approach to the limiting procedure for conditional measures and insightful discussions. We also thank Mikhail Sodin and Russell Lyons for helpful comments regarding the presentation.

\end{document}